\documentclass[a4paper]{amsart}
\usepackage{amscd, amssymb}
\usepackage{xypic}
\usepackage[matrix, arrow]{xy}
\usepackage{array}

\begin{document}



\def\rX{\mathrm{X}}
\def\coker{\operatorname{coker}}
\def\sign{\operatorname{sign}}
\def\O{\operatorname{O}}
\def\trace{\operatorname{trace}}
\def\oR{\operatorname{R}}
\def\GL{\operatorname{GL}}
\def\ind{\operatorname{ind}}
\def\ev{\operatorname{even}}
\def\odd{\operatorname{odd}}
\def\Fix{\operatorname{Fix}}
\def\Stab{\operatorname{Stab}}
\def\A{\mathcal{A}}
\def\H{\mathcal{H}}
\def\h{\mathcal{H}}
\def\E{\mathcal{E}}
\def\G{\mathcal{G}}
\def\F{\mathcal F}
\def\X{\mathcal X}
\def\mor{\operatorname{mor}}
\def\spin{\operatorname{spin}}
\def\hot{\hat{\otimes}}
\def\B{\mathbb B}
\def\cK{\mathbb K}

\newcommand*{\KK}{\mathrm{KK}}
\newcommand*{\RKK}{\mathrm{RKK}}
\newcommand*{\K}{\mathrm{K}}
\newcommand*{\Ktop}{\mathrm{K}^{\mathrm{top}}}
\newcommand*{\KX}{\mathrm{KX}}
\newcommand*{\RRKK}{\mathscr{R}\KK}

\def\RRKK{\mathcal{R}\KK}

\newcommand{\lam}{\Lambda}
\newcommand{\dulam}{\widehat{\Lambda}}

\newcommand{\Cc}{\mathcal{C}}

\newcommand*{\Comp}{{\mathbb{K}}}
\newcommand*{\Bound}{{\mathbb{B}}}
\newcommand*{\cross}{\mathbin{\rtimes}}

\newcommand{\T}{\mathbb{T}}
\newcommand{\Q}{\mathbb{Q}}
\newcommand{\Z}{\mathbb{Z}}
\newcommand{\C}{\mathbb{C}}
\newcommand{\N}{\mathbb{N}}
\newcommand{\R}{\mathbb{R}}
\newcommand{\Ktimes}{\times\!\!\!\!\!\times}

\newcommand{\Lambdah}{\widehat{\Lambda}}
\newcommand{\Deltah}{\widehat{\Delta}}
\newcommand{\ip}[1]{\langle #1 \rangle}

\newcommand{\Kn}{\mathcal{K}_{\nu}}
\newcommand{\Ko}{\mathcal{K}_{\omega}}
\newcommand{\Kob}{\mathcal{K}_{\bar{\omega}}}
\newcommand{\ltg}{\ell^2G}
\newcommand{\tG}{\tilde{G}}
\newcommand{\tN}{\tilde{N}}
\newcommand{\Go}{G_{\omega}}
\newcommand{\Co}{\mathbb{C}_{\omega}}
\newcommand{\ob}{\bar{\omega}}
\newcommand{\Gob}{G_{\bar{\omega}}}
\newcommand{\Cob}{\mathbb{C}_{\bar{\omega}}}
\newcommand{\ltgt}{\ltg_{\tau}}
\newcommand{\talpha}{\tilde{\alpha}}
\newcommand{\ttau}{\tilde{\tau}}
\newcommand{\FUpi}{\mathcal{F}_U\hot_{C_0(X),p_i}\Gamma_0(E)}
\newcommand{\FUpzero}{\mathcal{F}_U\hot_{C_0(X),p_0}\Gamma_0(E)}
\newcommand{\FUpone}{\mathcal{F}_U\hot_{C_0(X),p_1}\Gamma_0(E)}
\newcommand*{\defeq}{\mathrel{:=}}

\newcommand{\sour}{\mathscr{P}}
\newcommand{\PD}{\mathrm{PD}}
\newcommand{\pd}{\mathscr{PD}}
\newcommand{\lefschetz}{\mathscr{L}}
\newcommand{\Isom}{\mathrm{Isom}}
\newcommand{\Eul}{\mathrm{Eul}}
\newcommand{\Cliff}{\mathrm{Cliff}}
\newcommand*{\Ztwo}{{\mathbb Z/2}}
\newcommand{\dR}{\mathrm{dR}}
\newcommand{\comb}{\mathrm{comb}}
\newcommand{\Spinor}{\mathrm{Spinor}}
\newcommand*{\Dslash}{{\mathsf /\!\!\!\!D}}
\newcommand{\End}{\mathrm{End}}
\newcommand{\DD}{\mathscr{D}}
\newcommand*{\DslashS}{{\mathsf /\!\!\!\!D_{\mathrm{S}}}}
\newcommand*{\capslash}{{\; /\!\!\!\!\cap \;}}

\newcommand{\thmcite}[1]{{\bfseries\upshape \cite{#1}}}
\newcommand{\thmcitemore}[2]{{\bfseries\upshape \cite[#2]{#1}}}
\newcommand{\citemore}[2]{{\upshape \cite[#2]{#1}}}
\newcommand{\ltwoforms}{L^2(\Lambda^*_\C X)}
\newcommand{\ltwoformsg}{L^2(\Lambda^*_\C X)\hot \ell^2G}

\def\Cl{\mathcal{C}l}
\def\Aut{\operatorname{Aut}}
\def\Gx{X\rtimes G}
\def\Bott{\operatorname{Bott}}
\def\CC{\mathbb C}
\def\ClV{C_{\!{}_{V}}}
\def\comp{\operatorname{comp}}
\def\ddS{\stackrel{\scriptscriptstyle{o}}{S}}
\def\intW{\stackrel{\scriptscriptstyle{o}}{W}}
\def\intU{\stackrel{\scriptscriptstyle{o}}{U}}
\def\intT{\stackrel{\scriptscriptstyle{o}}{T}}
\def\E{\mathcal E}
\def\EE{\mathbb E}
\def\Egx{\mathcal{E}(G)\times X}
\def\EGG{\mathcal{E}(\mathcal{G})}
\def\EGT{\mathcal{E}(G/G_0)}
\def\EG{\underline{EG}}
\def\EGN{\mathcal{E}(G/N)}
\def\EH{\mathcal{E}(H)}
\def\I{{\operatorname{I}}}
\def\i{{\operatorname{i}}}
\def\Ind{\operatorname{Ind}}
\def\Id{\operatorname{Id}}
\def\infl{\operatorname{inf}}
\def\L{\mathcal L}
\def\lk{\langle}
\def\rk{\rangle}
\def\NN{\mathbb N}
\def\oplus{\bigoplus}
\def\pt{\operatorname{pt}}
\def\QQ{\mathbb Q}
\def\res{\operatorname{res}}
\def\sm{\backslash}
\def\top{\operatorname{top}}
\def\ZM{{\mathcal Z}M}
\def\ZZ{\mathbb Z}

\setcounter{section}{-1}

\def\Inf{\operatorname{Inf}}
\def\Ad{\operatorname{Ad}}
\def\id{\operatorname{id}}
\def\U{\mathcal U}
\def\PU{\mathcal PU}
\def\BG{\operatorname{BG}}
\def\eps{\epsilon}
\def\ev{\operatorname{ev}}
\def\om{\omega}
\def\Om{\Omega}
\def\F{\mathcal{F}}
\def\TT{\mathbb T}

\def\Flk{\F_{L,K}}
\def\mg{\mu_{G,A}}
\def\mgx{\mu_{\Gx,A}}
\def\ts{\hot}
\def\ga{\gamma}
\def\Br{\operatorname{Br}}
\def\Ab{\operatorname{Ab}}
\def\la{\lambda}
\def\gr{\operatorname{gr}}
\def\br{\Br^G_{\gr}(X)}
\def\mod{\operatorname{mod}}
\def\hato{\hat{\hot}}
\def\hotimes{{\hot}}
\def\hatoX{\hat{\hot}_{C_0(X)}}
\def\Ctau{C_{\tau}(X)}
\def\ctau{C_{\tau}}
\theoremstyle{plain}
 \newtheorem{theorem}{Theorem}[section]
\newtheorem{corollary}[theorem]{Corollary}
 \newtheorem{lemma}[theorem]{Lemma}
\newtheorem{lem}[theorem]{Lemma}
\newtheorem{prop}[theorem]{Proposition}
\newtheorem{proposition}[theorem]{Proposition} \newtheorem{lemdef}[theorem]{Lemma and Definition}
\theoremstyle{defn}
\newtheorem{defn}[theorem]{Definition}
\theoremstyle{definition}
\newtheorem{definition}[theorem]{Definition}
\newtheorem{defremark}[theorem]{\bf Definition and Remark}
\theoremstyle{remark}
\newtheorem{remark}[theorem]{Remark}
\newtheorem{assumption}[theorem]{Assumption}
\newtheorem{example}[theorem]{Example}
\newtheorem{note}[theorem]{Note}
\numberwithin{equation}{section} \emergencystretch 25pt
\renewcommand{\theenumi}{\roman{enumi}}
\renewcommand{\labelenumi}{(\theenumi)}

\newcommand*{\abs}[1]{\lvert#1\rvert}

\newcommand{\Rn}{\R^{n}}

\newcommand{\ltwo}{L^{2}}
\newcommand{\ltwotau}{L^2_\tau (F_\epsilon)}
\newcommand{\ltwotaugk}{L^2_\tau (F_{p,\epsilon,g}) }
\newcommand{\ltwotauk}{L^2_\tau (F_{p,\epsilon}) }
\newcommand{\ltwotaukv}{L^2_\tau (V_{p,\epsilon}) }
\newcommand{\ltwoformsrn}{\ltwo (\Lambda^{*}_\C \Rn)}

\newcommand{\repon}{\mathrm{R}\bigl(\mathrm{O}(n,\R) \bigr)}
\newcommand{\repgamma}{\mathrm{R}(\Gamma)}
\newcommand{\indaon}{\mathrm{ind}_{\mathrm{a}}^{\mathrm{O}(n,\R)}}
\newcommand{\indton}{\mathrm{ind}_{\mathrm{t}}^{\mathrm{O}(n,\R)}}
\newcommand{\indagamma}{\mathrm{ind}_{\mathrm{a}}^{\Gamma}}
\newcommand{\indtgamma}{\mathrm{ind}_{\mathrm{t}}^{\Gamma}}
\newcommand{\on}{\mathrm{O}(n,\R)}
\newcommand{\kgamma}{\K^{0}_\Gamma }
\newcommand{\ext}{\Lambda^{*}_\C}
\newcommand{\card}{\mathrm{card}}

\title[A Lefschetz fixed-point formula]
  {A Lefschetz fixed-point formula for certain orbifold C*-algebras}

\author[Echterhoff]{Siegfried
  Echterhoff}
   \address{Westf\"alische Wilhelms-Universit\"at M\"unster,
  Mathematisches Institut, Einsteinstr. 62 D-48149 M\"unster, Germany}
\email{echters@math.uni-muenster.de}

  \author[Emerson]{Heath
  Emerson}
  \address{Department of Mathematics and Statistics,
  University of Victoria,
  PO BOX 3045 STN CSCVictoria,
  B.C.Canada}
  \email{hemerson@math.uvic.ca}

  \author[Kim]{Hyun Jeong
  Kim}
 \address{Department of Mathematics and Statistics,
  University of Victoria,
  PO BOX 3045 STN CSCVictoria,
  B.C.Canada}
  \email{hjeong99@gmail.com}

\begin{abstract}
Using Poincar\'e duality in $\K$-theory, we state and prove a
Lefschetz fixed point formula for endomorphisms of crossed
product $C^{*}$-algebras $C_0(X)\cross G$ coming from covariant
pairs. Here $G$ is assumed countable, $X$ a manifold, and
$X\cross G$ cocompact and proper.
 The
formula in question expresses the graded
trace of the map on rationalized
$\K$-theory of $C_0(X)\cross G$ induced by the endomorphism, \emph{i.e.} the Lefschetz number,
in terms of fixed orbits and representation-theoretic data connected
with certain isotropy subgroups of the isotropy group at that point.

 \end{abstract}

\subjclass[2000]{19K35, 46L80}

\thanks{This research was supported by the EU-Network \emph{Quantum
  Spaces and Noncommutative Geometry} (Contract HPRN-CT-2002-00280)
  and the \emph{Deutsche Forschungsgemeinschaft} (SFB 478) and by the
  National Science and Engineering Research Council of Canada Discovery
  Grant program.}

\maketitle

\section{Introduction}

The goal of this article is to state and prove a
`noncommutative Lefschetz formula' for a certain class of
orbifold $C^{*}$-algebras $A$, and for a certain class of endomorphisms
$\alpha\colon A \to A$. The $C^{*}$-algebras in question are
the crossed products $A = C_0(X)\cross G$ where $X$ is a
manifold and $G$ is a countable group acting co-compactly and
properly on $X$. It is well-known that such actions give rise to 
orbifolds, and that the groupoids $X\rtimes G$ are Morita equivalent to 
the corresponding orbifold groupoids. 
The endomorphism $\alpha \colon A \to A$ is
 associated to a covariant pair $(\phi, \zeta)$ where
$\phi\colon X \to X$ is a map and $\zeta \in
\mathrm{Aut}(G)$ is a group automorphism, with $\zeta$ and
$\phi$ satisfying the equivariance condition $\phi \bigl(\zeta
(g)x\bigr) = g\phi (x)$. Note that this data corresponds to a self-map 
$\dot{\phi}\colon G\backslash X \to G\backslash X$ of 
the space of orbits, together with a coherent family of 
(finite) group homomorphisms, going between the isotropy 
groups attached to the orbits. It corresponds to an
 automorphism 
of the orbifold determined by the action of $G$ on $X$. 
 We consider the corresponding \emph{orbifold 
Lefschetz number}  taken by the trace of the induced map on the $\K$-theory
of the crossed-product $C^{*}$-algebra $A = C_0(X)\rtimes G$: 
\begin{equation}
\label{lefschetnumber}
\mathrm{Lef}(\alpha) \defeq \mathrm{trace}_s(\alpha_*\colon \K_*(A)_\Q \to \K_*(A)_\Q).
\end{equation}
The symbol $\mathrm{trace}_s$ denotes the graded trace (the trace on $\K_0$ minus the
trace on $\K_1$), and $\alpha_*$ is of course the map induced on $\K$-theory by the 
automorphism $\alpha\colon A \to A$.  We aim
to compute the Lefschetz number of $\alpha$ in \emph{geometric}
 terms. Specifically, we are going to compute it in terms of 1) the 
fixed orbits of the spatial map of the orbit space $G\backslash X$, and 
2) representation-theoretic information about the isotropy assigned to each 
such fixed orbit.

The geometry here
is therefore in some sense the geometry of the primitive ideal space of the crossed-product 
$C^{*}$-algebra $A$, which as a set
is a bundle over $G\backslash X$ with fibre over $\dot{x}\in G\backslash X$ the
irreducible dual of $\mathrm{Stab}_G(x)$, for any choice of $x \in \dot{x}$, but which
as a topological space has multiple points at orbits with nontrivial isotropy.

If
$G$ is trivial, or more generally, acts freely, then only fixed points of the 
induced map on the quotient $G\backslash X$ are detected, and we get 
essentially the classical Lefschetz fixed point theorem for $G\backslash X$. 

At the other extreme, where $X $ is trivial (is a point) and, hence $G$ is 
finite, we just have an automorphism of a finite group. Our Lefschetz theorem 
then recovers the following well-known fact about automorphisms of finite
groups: 
\begin{equation}
\#(\mathrm{Fix}(\hat{\zeta}\colon \widehat{G} \to \widehat{G} )\bigr) =
\frac{1}{|G|}\sum_{g \in G} |  \, Z_\zeta (g)|,
\end{equation}
where $Z_\zeta (g) = \{ h \in G \; | \; \zeta (h)g = gh\}$ and
$\hat{\zeta}\colon \widehat{G} \to \widehat{G} $ is the
permutation of the irreducible dual of $G$ induced by the
automorphism. This theorem is often expressed in terms of 
`twisted conjugacy classes' instead, the right hand side is 
trivially seen to be the number of $\zeta$-twisted conjugacy 
classes in $G$. 

In the general case, the local Lefschetz data around a fixed orbit 
in our theorem can be described as follows: the automorphism 
generates a family of subgroups of the isotropy group of each 
fixed point, and for each such subgroup, a virtual character of 
that subgroup. The characters are individually averaged, and the 
results added up.  

We now explain  this in a little more detail.

Let
 $\rho\colon \Gamma \to \mathrm{O}(n,\R)$ be an orthogonal representation of a finite
group $\Gamma$, and $A\in \mathrm{GL}(n,\R)$ a self-intertwiner of this representation; \emph{i.e.}
$A$ commutes with $\rho(\Gamma)$. Using this data we can define a conjugation-invariant
function
\begin{equation}
\label{eq:intro:orientation_character}\chi_{(\rho, A)} \colon \Gamma \to \Z, \;\;
\chi_{(\rho , A)} (g) = \mathrm{sign}\, \mathrm{det}(A_{|_{\mathrm{Fix}(g)}}),
\end{equation}
which we call the {\em orientation character} of the pair $(\rho,A)$.
Of course if $g\in \Gamma$ then $\mathrm{Fix}(g)$ is an $A$-invariant linear
subspace of $\R^{n}$ so this makes sense.

A pleasant and apparently not entirely 
obvious fact is that $\chi_{(\rho, A)}$ is a \emph{virtual character}, that is,
a difference of characters, of the group
$\Gamma$. We prove this. Indeed, this  `integrality result'
follows from index theory. It turns out that
 $\chi_{(\rho, A)}$ is the virtual
character associated
to the $\Gamma$-equivariant analytic index of the $\Gamma$-equivariant
Schrodinger-type operator obtained by perturbing the de Rham operator
$d+d^{*}$ on $L^{2}$-forms on $\R^{n}$ by the covector field $x\mapsto A\mathrm{X}$, where
$\mathrm{X}(x_1, \ldots , x_n) = x_1dx_1  +  \cdots + x_ndx_n$.

Of course, now the fact that $\chi_{(\rho, A)}$ is a character implies
that its average over the group with respect to normalized Haar measure is an integer, 
since by elementary representation theory this gives the dimension of the 
component of 
the trivial representation of the virtual representation corresponding to 
the virtual character $\chi_{(\rho, A)}$.

Returning to the general situation of $G$ acting on $X$, choose
a point $p$ from each fixed orbit of the induced map
$\dot{\phi}\colon G\backslash X \to G\backslash X$. For each
$p$ we have a secondary group action, and covariant pair, 
as follows. 

 Let ${L}_p \defeq \{ g\in G \; | \; \phi (gp) = p\}$;
then we have a group action of the isotropy group 
$\mathrm{Stab}_G(p)$ on ${L}_p$ by twisted
conjugation $h\cdot g \defeq \zeta (h)gh^{-1}$. 
Let the orbits of this action be represented by 
elements $g_1, \ldots , g_m$. For each $i$, let
 $\Gamma_{p,i}
\subset \mathrm{Stab}_G(p)$ be the stabilizer of $g_i$ under this
action.
 One easily checks that $\Gamma_{p,i}$ commutes
with $\phi \circ g_i$ and hence, differentiating and identifying
the tangent space at $p$ with $\R^{n}$, one obtains a
representation $\rho_{p,i}\colon \Gamma_{p,i}\to
\mathrm{O}(n,\R)$ and an intertwiner
 $A_{p,i} \defeq \mathrm{Id}- (\phi \circ g_i)'(p)$.
 Then our Lefschetz theorem reads as following:

\begin{theorem}
\label{intro:mainthm} In the above notation
\begin{equation}
\mathrm{Lef}([\alpha]) = \sum_{\dot{p}\in \mathrm{Fix}(\dot{\phi})}
\sum_{i }
\frac{1}{|\Gamma_{p,i} |   }
\sum_{h \in \Gamma_{p,i}}
\chi_{p,i} (h)
\end{equation}
where the $\chi_{\rho_{p,i}}$ are the index characters, as in \eqref{eq:intro:orientation_character}, 
so that 
$$\chi_{p, i}(h) = \mathrm{sign} \det (\mathrm{id} - D_p(\phi \circ g)_{|_{\mathrm{Fix}(h)}}).$$


\end{theorem}

The technique on which  the proof of our orbifold fixed point theorem relies on 
is quite
general, and can be phrased for general $C^{*}$-algebras: we
use the fact that for $C^{*}$-algebras satisfying the K\"unneth
theorem and the UCT and in addition satisfying Poincar\'e
duality in $\K$-theory, the Lefschetz number of an endomorphism
can be phrased as an index problem. This index problem 
arises from the automorphism and the cycles representing 
the fundamental classes of the duality. More precisely, the 
Lefschetz number can be realized as a Kasparov product in 
$\KK(\C, \C)$: one 
twists the fundamental class of the Poincar\'e duality by the automorphism, 
then pair with the dual fundamental class.  This index is computable 
in some situations by a local formula, as happens here. 
For more 
details of the general idea and an application to endomorphisms of Cuntz-Krieger 
algebras, see the preprint \cite{Emerson}.

That the $C^{*}$-algebras $C_0(X) \cross G$ and $\ctau
(X)\cross G$ are Poincar\'e dual is proved in \cite{EEK}. It
can be deduced from results of Kasparov on equivariant
$\KK$-theory. However, for purposes of applying the abstract
Lefschetz formula of \cite{Emerson} we need explicit
descriptions of the fundamental classes $\Delta$ and $\Deltah$.
The first part of the paper is devoted to finding such representatives.

In the second part, we analyse the orientation character
$\chi_{(\rho, A)}$ and in the third this becomes the critical
ingredient in the computation of the appropriate index pairing,
 which yields the Lefschetz theorem, Theorem \ref{intro:mainthm}.

The problem of giving a good description of the $\K$-theory of 
such orbifolds as appear here will be dealt with elsewhere. Such a 
description is needed to give a good formula for the 
\emph{global} Lefschetz number of an automorphism.  At the 
moment it seems to us to be a (surprisingly) delicate problem, 
however, and to keep down the length of the article, we have 
restricted our attention here to establishing the formula modulo 
$\K$-theory computations with a focus on the geometric, local 
description of our Lefschetz invariants.

\begin{note} All groups occurring in this paper are discrete. We
generally use group-algebra notation in connection with crossed
products. Thus, if $A$ is a $G$-$C^{*}$-algebra, then $A\cross
G$ is a completion of the linear span of the elements $a[g]$,
with $a \in A$ and $g\in G$.
\end{note}

\section{Fundamental classes}\label{fundamental}

Let $X$ be a complete Riemannian manifold and let $G$ be a
countable group acting isometrically, co-compactly and properly
on $X$. Let $\ctau (X)$ denote the algebra of continuous
sections of the Clifford algebra bundle of $X$. Since the group
$G$ acts isometrically on $X$, the action extends to an action
of $G$ on $\ctau (X)$. We can form the crossed product $\ctau
(X)\cross G$. To fix notation, we denote by 
$$d^x_g: T_{g^{-1}x}X\to T_xX$$
the differential of the action of $g$ on $X$ at the point $y=g^{-1}x$. 
It extends uniquely to a $*$-homomorphism 
$d^x_g: \Cl(T_{g^{-1}x}X)\to \Cl(T_xX)$  and the
action 
of $G$ on the Clifford bundle $C_\tau(X)$ is given by
 $$g(\varphi)(x)=d_g^x(\varphi(g^{-1}x)),$$
 for $\varphi\in C_\tau(X)$, $x\in X$ and $g\in G$. 
 
In this section, we are going to first review the
proof of the Poincar\'e duality between $C_0(X)\rtimes G$ and
$C_{\tau}(X)\rtimes G$, and then, using the proof, compute the
fundamental classes for this duality. Let us first recall the
following two equivalent definitions of Poincar\'e duality.

\begin{defn}[\emph{c.f.} \cite{EEK}, \cite{Em}]\label{duality} 
Let $\Lambda$ and $\Lambdah$ be $C^*$-algebras. Then $\Lambda$ and
$\Lambdah$ are \emph{Poincar\'e dual}
\begin{enumerate}
\item if there exist classes, called \emph{fundamental classes}, $\Delta \in \KK(\Lambda\hot \Lambdah,
\C)$ and $\Deltah \in \KK(\C, \Lambdah\hot \Lambda)$ such that
$\Deltah \hot_{\Lambdah} \Delta = 1_\Lambda$ and
$\Deltah\hot_\Lambda \Delta = 1_{\Lambdah}$, or equivalently,
\item if for every pair of $C^*$-algebras $A$ and $B$, there is an isomorphism
$$\Phi_{A,B}\,:\,\KK(\Lambda \hot A, B) \stackrel{\cong}{\longrightarrow}
\KK(A, \widehat{\Lambda} \hot B)$$ natural with respect to
intersection and composition products.
\end{enumerate}
\end{defn}

\begin{remark}\label{remark_definition}
It is easy to see the equivalence of the two definitions of
Poincar\'e duality. The isomorphism $\Phi_{A,B}$ of (ii) can be
obtained by the cap product with the class $\Deltah$ over $\Lambda$
and the inverse is given by the cap product with the class $\Delta$
over $\Lambdah$. On the other hand, for a given system of
isomorphisms $\{\Phi_{A,B}\}$, one can get classes
$\Delta=\Phi_{\Lambdah,\C}^{-1}(1_{\Lambdah})$ and
$\Deltah=\Phi_{\C,\Lambda}(1_{\Lambda})$.
\end{remark}

\begin{remark}\label{symmetry}
Note that when we say $\Lambda$ and $\widehat{\Lambda}$ are
\emph{Poincar\'e dual}, we already implicitly used the fact
that Poincar\'e duality is symmetric. Indeed one can show that
$\Delta':=\sigma_{}^*(\Delta)\in
\KK(\widehat{\Lambda}\hot\Lambda,\C)$ and
$\widehat{\Delta}':={\sigma_{}}_*(\widehat{\Delta})\in
\KK(\C,\Lambda\hot\widehat{\Lambda})$ satisfy Condition (i) in
Definition \ref{duality}, where $\sigma_{}$ is the flip
isomorphism.
\end{remark}

\begin{note}\label{dualnotation}
Under these circumstances, the maps
$$\Delta_*\colon \K_*(\Lambda) \to \K^*(\Lambdah), \;\;\;\; x \mapsto \hat{x}:=(x\otimes 1_{\Lambdah} ) \otimes_{\Lambda\otimes \Lambdah} \Delta$$ and
$$\Deltah_*\colon \K^*(\Lambdah) \to \K_*(\Lambda), \;\;\;\; y \mapsto \hat{y}:=\Deltah \otimes_{\Lambdah\otimes \Lambda} (y\hot1_\Lambda)$$
are inverse isomorphisms. Similarly, the maps
$$\Delta^*\colon \K_*(\Lambdah) \to \K^*(\Lambda), \;\;\;\; x \mapsto \hat{x}:=(1_\Lambda\otimes x ) \otimes _{\Lambda\otimes \Lambdah} \Delta $$ and
$$\Deltah^*\colon \K^*(\Lambda) \to \K_*(\Lambdah), \;\;\;\; y \mapsto \hat{y}:=\Deltah \otimes_{\Lambdah \otimes \Lambda} (1_{\Lambdah}\hot y)$$
are inverse isomorphisms.
\end{note}

Recall that Kasparov duality (see \cite{EmersonMeyer} for an
extensive discussion, or the original source
\citemore{Kas1}{Theorem 4.9}) states that, in this situation,
and more generally where $G$ is allowed to be locally compact,
there is a canonical isomorphism
\begin{equation}
\label{kasparovduality}K_{A,B}: \RKK^G(X; A,B) \stackrel{\cong}{\to}
\KK^G(\ctau (X)\hot A , B)
\end{equation} for any
$G$-$C^*$-algebras $A$ and $B$. If $G$ is discrete, then for $A$ and
$B$ equipped with the trivial $G$-action, we have the following
canonical isomorphism
\begin{equation}\label{covariant}
\begin{array}{rcl} C_{A,B}:
 \KK^G(C_{\tau}(X)\hot A, B)&\stackrel{\cong}{\to}&\KK\big((C_{\tau}(X)\rtimes G)\hot A, B\big);\\
\left[(\E,\varphi,F)\right]& \mapsto
&\left[(\E,\varphi\rtimes\pi,F)\right],
\end{array}
\end{equation}
where $\pi$ is the group representation on $\E$. Moreover, if
such $G$ acts properly on $X$, then, as a consequence of
\citemore{KasSk}{Theorem 5.4}, we have an isomorphism
\begin{equation}\label{descent}
E_{A,B}:\RKK^G(X; A,B)\stackrel{\cong}{\to} \KK\big(A,
(C_0(X)\rtimes G) \hot B\big).
\end{equation}

Combining all the isomorphisms above, we have Poincar\'e duality
between $C_0(X)\rtimes G$ and $C_{\tau}(X)\rtimes G$ as follows: for
all $C^*$-algebras $A$ and $B$ with trivial $G$-action, there exists
an isomorphism
\begin{equation}\label{poincare_duality}
\begin{array}{rcl}
\Phi_{A,B}:\KK\big((C_{\tau}(X)\rtimes G)\hot A, B\big)
&\stackrel{C_{A,B}^{-1}}{\longrightarrow}& \KK^G(C_{\tau}(X)\hot A, B)\\
&\stackrel{K_{A,B}^{-1}}{\longrightarrow}&\RKK^G(X; A,B)\\
&\stackrel{E_{A,B}}{\longrightarrow}&\KK\big(A, (C_0(X)\rtimes G)
\hot B\big),
\end{array}
\end{equation}
which is natural with respect to intersection and external products.

Now using the above system of isomorphisms $\{\Phi_{A,B}\}$ and
the equivalence of the two definitions of Poincar\'e duality
(see Definition \ref{duality} and Remark
\ref{remark_definition}) as well as the symmetry of Poincar\'e
duality (Remark \ref{symmetry}), we can compute fundamental
classes
\begin{equation}\label{fundamentalidentity}
\Delta=\sigma_{}^*\big(\Phi_{C_0(X)\rtimes G,\C}^{-1}(1_{C_0(X)\rtimes
G})\big)\,\,\,\,\,\,\text{ and }\,\,\,\,\,\,
\Deltah={\sigma_{}}_*\big(\Phi_{\C,C_{\tau}(X)\rtimes
G}(1_{C_{\tau}(X)\rtimes G})\big).
\end{equation}
For explicit descriptions for $\Delta$ and $\Deltah$, we need
an extensive discussion on the map $\Phi_{A,B}$, i.e., the maps
$C_{A,B}$, $K_{A,B}$ and $E_{A,B}$. We already know the map
$C_{A,B}$.
The map $K_{A,B}$ is the isomorphism of Kasparov's first
Poincar\'e duality. Recall the following Remark
\ref{remark_Kasparov} and Lemma \ref{firstpoincare} from
\cite{Kas1}.

\begin{remark}\label{remark_Kasparov}
\begin{enumerate}
\item Let $d:L^2(\Lambda_{\C}^*X)\rightarrow
    L^2(\Lambda_{\C}^*X)$ denote the (densely defined) de
    Rham operator. Let $D=d+d^*$ and let $F$ be the
    pseudodifferential operator
    ${D}(1+{D}^2)^{-\frac{1}{2}}$. Then
    $(L^2(\Lambda_{\C}^*X),F)$ is a cycle for an element in
    $\KK^{G}(\ctau (X), \C)$ where the action of
    $C_{\tau}(X)$ on $L^2(\Lambda_{\C}^*X)$ comes from the
    identification as vector bundles of the Clifford bundle
    of $X$ and the exterior bundle. We denote this cycle by
    $[D]$.

\item The map
$$\sigma_{X,C_\tau(X)}\colon \RKK^G(X; A,B) \to \KK^G(\ctau (X)\hot A , \ctau
(X)\hot B)$$ associates to a cycle $(\mathcal{E}, F)$ for
$\RKK^G(X;A,B)$ the cycle $(\ctau
(X)\hot_{C_0(X)}\mathcal{E}, 1\hot F)$ for $\KK^G(\ctau
(X)\hot A , \ctau (X)\hot B)$. The map $\sigma_{X,C_\tau(X)}$ is
natural with respect to intersection products in the sense
that
$$\sigma_{X,C_\tau(X)} (f\hot_B f') = \sigma_{X,C_\tau(X)}(f)\hot_{\ctau (X)\hot B}
\sigma_X(f'),$$
for $f\in \RKK^G(X; A,B), f'\in \RKK^G(X; B, C)$.

\item The map
$$
 p_X^*\colon \KK^G(A,B) \to \RKK^G(X; A,B),
$$ at the level of cycles, tensors with the standard
representative of $1_{C_0(X)} \in \KK^G(C_0(X), C_0(X))$.
Note that $p_X^*$ is natural with respect to intersection
products in the sense that $p_X^*(f_1\hot_B f_2) =
\big(p_X^*(f_1)\big)\hot_{X,B}\big(p_X^*(f_2)\big)$ for $f_1\in \KK^G(A, B)$,
$f_2\in \KK^G(B, C)$.

\item Let $\rho$ denote the metric on $X$. Then there
    exists an open neighbourhood $U$ of the diagonal in
    $X\times X$ where for every point $(x,y)\in U$ there
    exists a unique geodesic from $x$ to $y$. For such $U$ let $\F_U$ be
    the ideal of $C_0(X)\hot \ctau (X)$ of Clifford
    sections
supported on $U$. Since $G$ acts isometrically and
    cocompactly on $X$, there exists $\epsilon>0$ such that
    $U_{\epsilon}:=\{(x,y)|\rho(x,y)<\epsilon\}$ is
    contained in the set $U$. Let
    $\theta_{\epsilon}(x,y)=\frac{\rho(x,y)}{\epsilon}(d_y\rho)(x,y)$.
     Then $(\F_{U_{\epsilon}},\theta_{\epsilon})$ defines a cycle in $\RKK^G(X; \C,
     \ctau (X))$ with $\theta_{\epsilon}$ as a multiplicative operator
     and $[(\F_{U_{\epsilon}},\theta_{\epsilon})]=[(\F_{U_{\epsilon'}},\theta_{\epsilon'})]$
    for any $0<\epsilon'\leq\epsilon$.  We denote the class
    $[(\F_{U_{\epsilon}},\theta_{\epsilon})]$ by $\Theta$, and we shall
    simply write $(\F_U,\theta)$ if we do not  want to specify the special 
    $\epsilon$ used in the construction.
\end{enumerate}

\end{remark}

The following is a special case of Kasparov's  \cite[Theorem 4.9]{Kas1}:

\begin{lem}\label{firstpoincare} Let $G$ act isometrically and cocompactly
 on a complete Riemannian manifold $X$. The
composition
{\small \begin{equation}
  \label{pdinverse}
 K_{A,B}:\RKK^G(X; A,B) \stackrel{\sigma_{X,C_\tau(X)}}{\longrightarrow}
 \KK^G(\ctau(X)\hot A , \ctau (X)\hot B)
 \stackrel{-\otimes [D]}{\longrightarrow} \KK^G(\ctau (X)\hot A , B)
\end{equation}}
is an isomorphism with inverse the composition
\begin{equation}
  \label{pd}
  K_{A,B}^{-1}:\KK^G(\ctau (X)\hot A , B) \stackrel{p_X^*}{\longrightarrow}
  \RKK^G(X; \ctau (X)\hot A , B)
  \stackrel{\Theta \hot - }{\longrightarrow}
  \RKK^G(X; A, B).
\end{equation}
\end{lem}

The map $E_{A,B}$ is the isomorphism from \cite[Theorem 5.4]{KasSk}. To
understand the map $E_{A,B}$ explicitly, we need to understand
two ingredients. Firstly, the \emph{descent homomorphism}
\begin{equation}\label{descent}
 \mathrm{descent}\colon \RKK^G(X; A,B) \to \KK(C_0(X,A)\cross G,
 C_0(X,B)\cross G).
 \end{equation}
Secondly, the \emph{unit class} $[E] \in \K_0( C_0(X)\cross G)$,
defined whenever $G\backslash X$ is compact: $[E]$ is defined by the
finitely generated projective right $C_0(X)\cross G$-module $E$ which is
the completion of $C_c(X)$ with respect to the inner product
\begin{equation}\label{innerproductonE}
 {\langle}\varphi, \varphi'{\rangle}(x,g) = \varphi (x)\varphi'(gx).
\end{equation}
For future reference, the right action of $C_0(X)\cross G$ on
 $E$ is given by
\begin{equation}
 \varphi f \, (x) = \varphi (x)f(x), \;\;\;\; \varphi g\, (x) = \varphi (gx), \;\;\;\;\;
 g\in G, f\in C_0(X).
\end{equation}

\begin{remark} $[E]$ is also represented by the projection
 $P \in C_0(X)\cross G$, $$P = \sum_{g \in G} \varphi g(\varphi),$$ where
 $\varphi \in C_c(X)$ is chosen so that
  $0\le \varphi \le 1$,
 $\sum_{g \in G} g(\varphi)^2 = 1.$ See \cite{BCH}.
\end{remark}

\begin{lem}[\citemore{KasSk}{Theorem 5.4}] Let $G$ act on $X$ isometrically and properly. The map
$E_{A,B}:
 \RKK^G(X; A,B) \to \KK(A,(C_0(X)\cross G)\hot B)$ given by the composition
\begin{multline}\label{descentisomorphism}
 \RKK^G(X; A, B) \stackrel{\mathrm{descent}}{\longrightarrow}
 \KK(C_0(X,A)\cross G, C_0(X,B)\cross G) \\ \stackrel{[E]\otimes - }
 {\longrightarrow} \KK(A, (C_0(X)\cross G)\hot B)
\end{multline}
is an isomorphism whenever $A$ and $B$ are $G$-trivial
$C^*$-algebras.
\end{lem}

\begin{remark}[c.f. \cite{KasSk}]\label{explicitdescriptionofdescentisomorphism}
The map $E_{A,B}$ can be explicitly described as follows.
Suppose we have a cycle $(\mathcal{E}, F)$ for $\RKK^G(X;
A,B)$. Then $\mathcal{E}$ is a right $C_0(X,B)$-module, and a
left $C_0(X,A)$-module, and the two actions of $C_0(X)$ on the
left and right agree. Furthermore, the group  $G$ acts on
$\mathcal{E}$. We can assume by averaging that $F$ is exactly
$G$-invariant. Now we complete the compactly supported elements
of $\mathcal{E}$ to a right $C_0(X,B)\cross G$-module
$\tilde{\mathcal{E}}$ using the inner product valued in
$C_0(X,B)\cross G$,
\begin{equation}
 {\langle}\xi,\xi'{\rangle} =\sum_{h\in G} {\langle}\xi, h(\xi'){\rangle}[h].
\end{equation}
The right action of $C_0(X,A)\cross G$ is given by letting
$C_0(X,A)$ act as originally, and $G$ acting by $\xi h =
h^{-1}(\xi).$

Finally, we note that \emph{adjointable operators on the right
$C_0(X,B)\cross G$-module $\tilde{\mathcal{E}}$ are in
$1$-$1$-correspondence with $G$-equivariant operators on $\E$. }
(Generalizations of the isomorphism $E_{A,B}$ are given in \cite{EM, EEK}.)
\end{remark}

Now we are going to construct fundamental classes $\Delta$ and
$\Deltah$.
First, we construct the class $\Delta$. Recall Remark
\ref{remark_Kasparov} (i) for the discussion below. We define two
commuting $*$-homomorphisms $\ctau(X)\cross G \to
\Bound(\ltwoformsg)$ and $C_0(X)\cross G \to \Bound (\ltwoformsg)$,
by constructing two corresponding covariant pairs. We let $e_g \in
\ell^2G$ denote the point-mass at $g \in G$.

The $*$-homomorphism $\ctau (X)\cross G \to \Bound(\ltwoformsg)$ is
defined by the covariant pair
\begin{equation}\label{representationCtau}\varphi(\xi \hot e_g) = \varphi \cdot
\xi \hot e_g, \;\;\;\;\; h(\xi\hot e_g) = h\cdot \xi \hot
e_{gh^{-1}}, \end{equation} for $g,h \in G, \xi \in \ltwoforms,
\varphi \in \ctau (X)$. The $*$-homomorphism $C_0(X)\cross G \to
\Bound(\ltwoformsg)$ is defined by the covariant pair
\begin{equation}\label{representationCzero}
f(\xi \hot e_g) = g^{-1}(f)\cdot \xi \hot e_g, \;\;\;\; h(\xi\hot
e_g) = \xi\hot e_{hg},\end{equation}
for $ f\in C_0(X), g,h\in G,
\xi \in \ltwoforms$. The dots indicate the actions already implied
in the cycle $D=(\ltwoforms , F)$ of Kasparov (see Remark
\ref{remark_Kasparov} (i)); note that $C_0(X)$ embeds in $\ctau
(X)$. Observe that the two $*$-homomorphisms just defined
 commute, and so determine a $*$-homomorphism
$$\Pi: C_0(X)\cross G\hot \ctau (X)\cross G \to \Bound(\ltwoformsg).$$

Let $g\in G$, $f\in C_0(X)$ and $\varphi \in \ctau (X)$ be compactly
supported. If $T$ is a locally compact operator on $\ltwoforms$, \emph{e.g.} if
$T = F^{2} - 1$, then
$$(f\hot\varphi)(T\hot 1)(\xi\hot e_g) = g^{-1}(f)\varphi T\xi \hot e_g,$$ so that
$(f\hot\varphi)(T\hot 1)$ acts as the block diagonal operator
$\oplus_{g\in G} g^{-1}(f)\varphi (T\hot1)$, which has compact
blocks. As
$$g^{-1}(f)\varphi = 0 \;\;\;\text{for}\; g\notin H \defeq \{h \in G \mid h^{-1}\bigl(
\mathrm{supp}(f)\bigr)\cap \mathrm{supp} (\varphi)\not= \emptyset\},
$$
and since the indicated set $H$ is finite, since $G$ acts properly,
there are only finitely many blocks. Thus $(f\hot\varphi)\bigl(
(F^2-1)\hot 1\bigr)$ is compact. This observation and similar ones
prove that the Hilbert space $\ltwoformsg$ equipped with the
$*$-homomorphism
$\Pi:\ctau (X)\cross G \otimes C_0(X)\cross G \to \Bound(\ltwoformsg)$
defined above and the operator $F\hot 1$ defines a
cycle for $\KK(C_0(X)\cross G\hot\ctau (X)\cross G , \C)$.

\begin{defn}\label{classdelta}
We define $$\Delta \in \KK(C_0(X)\cross G \hot\ctau (X)\cross
G, \C)$$ to be the class of the cycle $(\ltwoformsg,\,\Pi,\, F\hot
1)$ above.
\end{defn}

We shall show below that $\Delta$ coincides with  the fundamental class
of (\ref{fundamentalidentity}).
\\

The dual class $\Deltah$ is more complicated to write down.
Recall the class $\Theta=[(\F_U,\theta)]$ from Remark
\ref{remark_Kasparov} (iv). We consider the completion $\E$ of
$C_c(X){\F_U}\hot \C G$ (with action of $C_c(X)$ with respect to the first variable in $U\subseteq X\times X$) equipped with the following $ \ctau
(X)\cross G\hot C_0(X)\cross G $-valued inner product:
\begin{equation}\label{innerproduct}{\langle}\alpha \hot [g], \alpha'\hot [g']{\rangle} = \sum_{h
\in G} g^{-1}\bigl(\alpha^{*}h(\alpha ')\bigr) [g^{-1}hg']\hot [g^{-1}h].
\end{equation}
Here $[g^{-1}h]$ is understood as in $G \subset C_0(X)\cross G$,
$[g^{-1}hg']$ is understood as in $G \subset \ctau (X)\cross G$
and $g^{-1}\bigl(\alpha^{*}h(\alpha ')\bigr)$ is understood as
in $C_{\tau}(X)\hot C_0(X)$, equipped with the diagonal $G$-action. 
The right module-structure is
given by
\begin{equation}\label{rightaction}\begin{array}{rcl} 
&&(\alpha \hot [g])f = f\alpha
\hot [g],\quad\quad(\alpha \hot [g])h = h^{-1}(\alpha) \hot
[h^{-1}g],\\
&&(\alpha \hot [g])\varphi = \alpha g(\varphi) \hot [g],\;\; (\alpha
\hot [g])h' = \alpha \hot [gh'],
\end{array}
\end{equation}
where $f\in C_0(X), h \in G\subset C_0(X)\cross G$ and $\varphi
\in \ctau (X), h' \in G \subset \ctau (X)\cross G$.

Note that any $G$-invariant element of $M(\ctau (X)\hot C_0(X))$ acts
as an operator on $\E$ by multiplication in the ${\F_U}$-variable.
The $G$-invariance is needed to commute with the action of $G\subset
C_0(X)\cross G$ on the right.

\begin{defn}
\label{classdeltah} The class $\Deltah \in \KK(\C, \ctau (X)\cross
G\hot C_0(X)\cross G)$ is given by the pair $(\E, \theta_G)$, where
we let $\theta_G$ be the operator on $\E$ induced by the
$G$-invariant multiplier $\theta$ of ${\F_U} \subset C_0(X)\hot
\ctau (X)$ described above, and $\E$ is the Hilbert module as above
with the inner product and the right actions given in
\eqref{innerproduct} and \eqref{rightaction}.
\end{defn}

Now we need to prove that the classes $\Delta$ and $\Deltah$ defined
above are actually the fundamental classes for
$\Lambda=C_0(X)\rtimes G$ and $\Lambdah=C_{\tau}(X)\rtimes G$, i.e.,
they satisfy identities $\Delta=\sigma_{}^*\big(\Phi_{C_0(X)\rtimes
G,\C}^{-1}(1_{C_0(X)\rtimes G})\big)\in \KK(\Lambda\hot\Lambdah,\C)$ and
$\Deltah={\sigma_{}}_*\big(\Phi_{\C,C_{\tau}(X)\rtimes
G}(1_{C_{\tau}(X)\rtimes G})\big)\in \KK(\C,\Lambdah\hot\Lambda)\big)$ from
(\ref{fundamentalidentity}), where $\sigma$ denotes the flip isomorphism.
This will follow from

\begin{prop}\label{imagesofone}
Let $\Delta$ and $\Deltah$ be the classes defined in Definitions
\ref{classdelta} and \ref{classdeltah} and let $\Phi_{(-, -)}$ be
the isomorphisms \eqref{poincare_duality}. Then
\begin{equation}
  {\sigma_{}}_*\big(\Phi_{\C, \ctau (X)\cross G} (1_{\ctau (X)\cross G})\big) = \Deltah, \;\;
  \text{and } \;\;  \Phi_{C_0(X)\cross G, \C} (\sigma_{}^*\Delta) = 1_{C_0(X)\cross G}.
  \end{equation}
\end{prop}
For the proof, we need some preliminary discussion.

\begin{defremark}\label{multiplicationmap}
We give -- here and elsewhere -- the crossed product $C_0(X)\cross
G$ the structure of a trivial $G$-$C^*$-algebra. Let $A$ be an
$X\cross G$-algebra. The \emph{multiplication class} $m_A$ is the
class
$$m_A\in \RRKK^G(X; A\hot C_0(X)\cross G , A)$$
given by the cycle $(A\hot \ell^2G, 0)$, where the right Hilbert
$A$-structure of $A\hot \ell^{2}G$ is the obvious one, and the
further module structures are as follows.

Note that there are two $G$'s involved here; one the $G$ which
appears in $\RRKK^{G}$, the other which appears in the crossed
product $C_0(X)\cross G$. To reduce confusion, we refer to the
action of the former as the \emph{equivariant} action. The
equivariant action of $G$ on $A\hot \ell^{2}G$ is then given by
$h(a\hot e_g) = h(a)\hot e_{gh^{-1}}$. The $C_0(X)$-structure
is by multiplication in the $A$ factor.

The representation of the crossed product is given by the covariant
pair
$$f(a\hot e_g) = g^{-1}(f) a\hot e_g,\;\;\;\; h(a\hot e_g) = a\hot e_{hg}.
$$
One easily checks that this is a covariant pair. The left actions of
$C_0(X)\cross G$ and of $C_0(X)$ clearly commute, and, finally, the
$*$-homomorphism $C_0(X)\cross G \to
\Bound\bigl( A\hot \ell^{2}G)$ is equivariant -- that is
$C_0(X)\cross G$ acts as $G$-invariant operators on
$A\hot \ell^{2}G$.
Since the action of $A\hot C_0(X)\rtimes G$ is by compact operators on 
$A\hot \ell^2G$ we get a cycle as required.
\end{defremark}

\begin{lem}\label{mA} The isomorphism ${E_{C_0(X)\rtimes G,\C}}\colon
 \RKK^G(X; C_0(X)\cross G, \C) \to \KK(C_0(X)\cross G , C_0(X)\cross G)$
 maps $m_{C_0(X)}$ to $1_{C_0(X)\cross G}$.
\end{lem}

\begin{proof}
Apply the explicit description in Remark \ref{explicitdescriptionofdescentisomorphism}.
We get the
completion of $C_0(X)\hot \C G$ with respect to the
following inner product:
$$
\langle a\hot e_g, a'\hot e_{g'}\rangle = a^{*}g^{-1}\bigl(g'
(a')\bigr) [g^{-1}g'] \in C_0(X)\cross G, \;\;\;\; \text{for
}a,a'\in C_0(X), g,g'\in G
$$
and the right $C_0(X)\cross G$-module structure
$$
(a\hot e_g)f = fa\hot e_g,\;\;\;\;\;\;\;\; (a\hot e_g)h =
h^{-1}(a)\hot e_{gh}.$$ The left action of $C_0(X)\cross G$ is given
by
$$
f(a\hot e_g) = g^{-1}(f)a\hot e_g,\;\;\;\;\; h(a\hot e_g) = a\hot e_{hg}.
$$

Let $\overline{ C_0(X)\hot \ell^{2}G}$ be the completion of the
above to a Hilbert module. We define a map $W\colon \overline{
C_0(X)\hot \ell^{2}G} \to C_0(X)\cross G$, where the co-domain has
its standard $C_0(X)\cross G$-bimodule structure, by the formula
$$W(a\hot e_g) = g(a)[g].$$
Then

$$\langle W^{-1}(a[g]), W^{-1}(a'[g'])\rangle = \langle g^{-1}(a)\hot e_g, (g')^{-1}(a') \hot e_{g'}\rangle
= g^{-1}(a^{*}a') [g^{-1}g'],$$ and
$$\begin{array}{c}
W^{-1}(a[g]h) = W^{-1}(a[gh]) = h^{-1}\bigl( g^{-1}(a)\bigr)\hot e_{gh} = \bigl( g^{-1}(a)\hot e_g\bigr)h, \\
W^{-1}(a[g]f) = W^{-1}(ag(f)[g]) = g^{-1}(a)f \hot e_g = \bigl(
W^{-1}(a[g])\bigr)f.
\end{array}$$

Hence $W$ gives an isometry between the inner product we have
defined initially, on $C_0(X)\hot \C G$, and the usual inner product
on the free, rank one Hilbert $C_0(X)\cross G$-module.

Similarly, one checks that
$W$ conjugates the left $C_0(X)\cross G$-module structure we have
defined above, and the standard one by algebra multiplication.

Therefore the image of the class $m_{C_0(X)}$ under the map
${E_{C_0(X)\rtimes G,\C}}$ sends the cycle for the multiplication
class, to a cycle which is unitarily equivalent to the standard
representative of $1_{C_0(X)}$, so that ${E_{C_0(X)\rtimes G,\C}}
(m_{C_0(X)}) = 1_{C_0(X)}$ as claimed.
\end{proof}

\begin{proof}[Proof of Proposition \ref{imagesofone}] Consider first
the fundamental class $\Delta$. This is accomplished by a direct
computation:
by Lemma \ref{mA}, it remains to apply the map, from
Remark \ref{remark_Kasparov} (ii),
$$\sigma_{X,\ctau(X)}\colon \RKK^{G}(X; C_0(X)\cross G , \C) \to \KK^{G}(\ctau (X)\hot C_0(X)\cross G , \ctau (X))$$
to the class $m_{C_0(X)}$. A straightforward application of the definition gives the
cycle $\bigl( \ctau (X)\hot \ell^{2}G, 0\bigr)$. The left action of
$\ctau (X)$ is given by
$\varphi (a\hot e_{g}) = \varphi a\hot e_{g},$
the group $G$ acts by
$h(a\hot e_{g}) = h(a)\hot e_{gh^{-1}}.$ The left action of $C_0(X)\cross G$ is given by the
covariant pair
$$f(a\hot e_g) = g^{-1}(f)a\hot e_g, \;\;\;\; h(a\hot e_g) = a \hot e_{hg}.$$

Finally, we take the product of the class of this cycle, with
the class $D \in \KK^{G}(\ctau (X), \C)$ of Kasparov (see
Remark \ref{remark_Kasparov} (i)). Comparing to Definition
\ref{classdelta}, we see that the modules are the same. The
axioms for a Kasparov product imply that the operator involved
in the product is also that described in Definition
\ref{classdelta}.

The assertion regarding $\Deltah$ is similar but slightly more straightforward; we leave
its confirmation to the reader.
\end{proof}

\section{The orientation character}\label{schrodinger}

Before proceeding to the Lefschetz theorem, we need to perform
an index calculation on Euclidean space $\Rn$ generalising the
computation of the index of the Schrodinger operator
$\frac{d}{dx} \pm x$ on $\ltwo(\R)$ (see \cite{Higson-Roe}). 
The analogue for $\Rn$ of the harmonic oscillator in dimension $1$
is the operator $D+ \rX$, where $D = d+ d^{*}$ is the de Rham
operator acting on $\ltwoformsrn$, and $\rX$ is Clifford
multiplication by the form $ x_1dx_1 + \cdots + x_ndx_n$ on $\Rn$.
The (unbounded) cycle $(\ltwoforms, D+\rX)$ represents the Kasparov product
$$[\mathrm{X}]\hot_{\ctau (\Rn)}[D]  \in \KK^{\mathrm{O}(n,\R)}(\C, \C) = \repon$$
of the class $[D]$ of the de Rham operator (see Remark \ref{remark_Kasparov})
the generator $[\mathrm{X}]\in
\KK^{\mathrm{O}(n,\R)}(\C, \ctau (\Rn))$ constructed  via the vector field $\rX$ as 
explained below. 
It is the content of Kasparov's Bott-periodicity theorem 
\cite[Theorem 7 of \S 5]{Kas-operator}
that $[\mathrm{X}]\hot_{\ctau (\Rn)}[D]=1\in \repon$.

In this section we extend Kasparov's calculations to the following 
more general situation: we assume that $\Gamma$ is a compact group acting 
on $\R^n$ via an orthogonal representation $\rho: \Gamma\to \on$.
Moreover, we shall assume that $A\in \GL(n,\R)$ commutes with $\rho$. 
We then obtain a $\Gamma$-invariant 
Fredholm operator $A\rX+D$ (or a bounded version of it)
and we need to compute the $\Gamma$-index 
$$\operatorname{index}^\Gamma(A\rX+D)\in \repgamma.$$
We shall do this in two different ways:
in a first version we make use of Kasparov's ideas for the proof of his Bott-periodicity theorem 
by reducing the computations to appropriate two- and one-dimensional subspaces.
In a second version we sketch the argument how the result can also be obtained 
from  a use of the Atiyah-Singer Index Theorem for open subsets of $\R^N$
(see \cite{AS1}) together with some calculations  
given by Atiyah and Segal in \cite{AS2}.

Before we do this we need to recall
the relation between vector fields on a manifold $X$ 
and corresponding classes in $K_0^\Gamma(C_\tau(X))$.
So let $X$ be any Riemannian manifold such that the compact group
$\Gamma$ acts isometrically on $X$. Suppose that
$v:X\to TX$ is a $\Gamma$-invariant continuous vector field on $X$ 
such that there exists a compact set $K\subseteq X$ with $v(x)\neq 0$ 
outside $K$. We then say that $v$ is {\em co-compactly supported}.
From $v$ we construct a new vector field $\tilde{v}$ as follows:
choose a $\Gamma$-invariant positive continuous function $\varphi:X\to [0,1]$ such that 
$\varphi\equiv 0$ on $K$ and such that $1-\varphi\in C_0(X)$. Then set
$$\tilde{v}(x):=\varphi(x)\frac{v(x)}{\|v(x)\|}.$$
The vector field $\tilde{v}$ 
acts as self-adjoint bounded operator
 on $C_{\tau}(X)$ by point-wise Clifford multiplication 
such that $\tilde{v}^2-1$ (which is point-wise multiplication by $x\mapsto( \|\tilde{v}(x)\|^2-1)$)
lies in $C_\tau(X)$. It thus defines a class 
$[v]\in \KK^\Gamma(\C, C_\tau(X))=K_0^\Gamma(C_\tau(X))$.

Two such vector fields $v_0,v_1: X\to TX$ are said to be homotopic, if there 
exists a co-compactly supported 
$\Gamma$-invariant continuous map $v: X\times [0,1]\to TX$  such that
$v(x,t)\in T_xX$ for all $(x,t)\in X\times[0,1]$ and 
 $$v|_{X\times \{0\}}=v_0\quad\text{and}\quad v|_{X\times\{1\}}=v_1.$$

\begin{lemma}\label{lem-homotopy}
Suppose that $v:X\to TX$ is a co-compactly supported $\Gamma$-invariant 
continuous vector field. 
Then the class $[v]\in \K_0^\Gamma(C_\tau(X))$ does not depend on the choice 
of the function $\varphi$. Moreover, two homotopic co-compactly supported 
vector fields on $X$ determine the same class in $\K_0^{\Gamma}(C_\tau(X))$.
\end{lemma}
\begin{proof} Suppose that $v(x)\neq 0$ outside the compact set $K\subseteq X$
and suppose that 
$\varphi_0$ and $\varphi_1$ are two functions which vanish on $K$ and which 
have value $1$ at $\infty$..
Then 
$$t\mapsto \tilde{v}_t=(t\varphi_1+(1-t)\varphi_0)\frac{v}{\|v\|}$$
 is an operator homotopy 
between $\tilde{v}_0$ and $\tilde{v}_1$ which proves the first assertion. 
A similar argument gives the second assertion.
\end{proof}

Recall from  Remark \ref{remark_Kasparov}
the construction of the Dirac class $[D]=[D_X]\in \KK^\Gamma_0(C_\tau(X),\C)$ 
given by the de Rham operator $D=d+d^*: L^2(\Lambda_\C^*(X))\to L^2(\Lambda_\C^*(X))$.
Note that if $U\subseteq X$ is any open $\Gamma$-invariant
sub-manifold, then $[D_X]$ restricts to the 
class $[D_U]$ under the canonical inclusion $\iota_U:C_{\tau}(U)\to C_{\tau}(X)$.
The following basic (and certainly well-known)
 lemma turns out to be extremely useful for our computations.

\begin{lemma}\label{lem-support}
Suppose that $v:X\to TX$ is a co-compactly supported $\Gamma$-invariant
vector field on $X$. Let $K\subseteq X$ be compact such that $v$ does not vanish 
outside $K$ and let $U\subseteq X$ be an open $\Gamma$-invariant neighborhood of 
$K$ in $X$. Then 
$$[v_U]\otimes_{C_\tau(U)}[D_U]=[v]\otimes_{C_\tau(X)}[D]\in \repgamma,$$
where $v_U:U\to TU$ denotes the restriction of $v$ to $U$.
\end{lemma}
\begin{proof}
By the construction of the class $[v]$ we may assume without loss of generality that
there exists a compact $\Gamma$-invariant set $C\subseteq U$ such that $\|v(x)\|=1$ 
for all $x\notin C$.  The Kasparov product
$[v]\otimes_{C_\tau(X)}[D]$ is represented by the pair
\begin{equation}\label{Kasp-product}
(L^2( \Lambda^*_\C(X)), T)
\quad\text{with}\; T=\lambda_{v(x)}+\lambda_{v(x)}^*+\sqrt{1-\|v(x)\|^2}\frac{D}{\sqrt{1+D^2}}
\end{equation}
with $D=d+d^*$,
which can be deduced from \cite[Remark 3 on p. 541]{Kas-operator}.
Since $\|v(x)\|^2=1$ outside $U$, it follows that the second summand vanishes
on $X\setminus U$. 
It is then clear that
 $L^2(\Lambda^*_\C(X))$ decomposes into the direct product of $T$-invariant subspaces
$L^2(\Lambda^*_\C(U))\oplus L^2(\Lambda^*_\C(X\setminus U))$ such that the restriction of $T$ to $L^2(\Lambda^*_\C(U))$
gives the product $[v_U]\otimes_{C_\tau(U)}[D_U]$ (since $D$ is local). The restriction 
of $T$ to $L^2(\Lambda^*_\C(X\setminus U))$ is given  point-wise by the unitary operator
$\lambda_{v(x)}+\lambda_{v(x)}^*$ (it is unitary since $\|v(x)\|=1$) and hence has index $0$.
\end{proof}

\begin{remark}\label{rem-product} 
 Suppose that $(\E_1, \phi_1, F_1)$ and $(\E_2,\phi_2, F_2)$ are  two Kasparov cycles
 giving elements $x\in  \KK^G(A,B)$ and $y\in \KK^G(B, C)$, respectively,
 where we assume here that  $G$ is a compact group.
 Assume that both operators $F_1, F_2$ are $G$-invariant and self-adjoint with $\|F_1\|\leq 1$.
Suppose further that  $F\in \B(E_1\hot_B\E_2)$ is a self-adjoint  $F_2$-connection, i.e.,
$$\Theta_\xi F_2 - (-1)^{\deg(\xi)\cdot\deg(F_2)}F\Theta_\xi\in \cK(\E_2, \E_1\hot_B\E_2)$$
for all $\xi\in \E_1$, where $\Theta_{\xi}:\E_2\to \E_1\hot_B\E_2; \eta\mapsto \xi\hot_B\eta$. 
Let 
$$T=(F_1\otimes 1)+\sqrt{1- F_1^2\otimes 1}^{1/2}F\in \B(\E_1\hot_B\E_2).$$
It follows then from \cite[18.10.1]{Blackadar} that $(\E_1\otimes_B\E_2, \phi_1\otimes 1, T)$
is a representative for the Kasparov product $x\hot_By\in \KK^G(A,C)$, provided 
$[T,\phi_1(A)\hot 1]\in \cK(\E_1\hot_B\E_2)$. Formula (\ref{Kasp-product}) is a direct consequence
of this principle. But we shall use this principle also  in a more advanced setting in \S3 below.
\end{remark}

We now  specialize to the case where $X=V$ is a finite dimensional Euclidean vector space
together with a linear action 
$\rho:\Gamma\to \O(V)$. We want to
give explicit computations of the product $[v]\otimes_{C_\tau(V)}[D_V]$
in case where $v:V\to TV=V\times V$ is given by 
$v(x)=Ax$ for some $A\in \operatorname{GL}(V)$ which commutes with the representation $\rho$.
We shall always write $A\rX$ for this vector field. We shall 
show below that the product 
$[A\rX]\otimes_{C_{\tau}(V)}[D]\in \KK^\Gamma(\C,\C)=\repgamma$
is equal to the orientation character $\chi_{(\rho,A)}$ as in

\begin{defn}\label{chi}
Let $\rho\colon \Gamma \to \O(V)$ and  $A \in \mathrm{GL}(V)$ as above.
The \emph{orientation character}
$\chi_{(\rho, A)}\colon \Gamma \to \Z$ is the conjugation-invariant function on $\Gamma $
$$\chi_{(\rho, A)} (g) \defeq \mathrm{sign}\, \mathrm{det} (A_{|_{\mathrm{Fix}(g)}}),$$
where $\mathrm{Fix}(g)\subseteq V$ denotes the space of fixed-points  for $g\in \Gamma$.
\end{defn}

The  set $\mathrm{Fix}(g)$  is of course a linear subspace of $V$ invariant under
$A$, so the formula makes sense. It is clear that $\chi_{(\rho,A)}$
is conjugation-invariant. 
The remaining part of this section is devoted to the proof of

\begin{theorem}
\label{charactertheorem} The orientation character $\chi_{(\rho, A)}$ is a virtual character of $\Gamma$ (i.e., a difference of
two characters). Under the identification of $\repgamma$ as the ring of $\Z$-linear combinations of
characters of $\Gamma$, we have
$$\chi_{(\rho,A)}=[A\rX]\otimes_{C_{\tau}(V)}[D].$$
\end{theorem}

\begin{remark}\label{rem-character}
(a)  Recall that the identification of $\repgamma$ with 
the ring of $\Z$-linear combinations of
characters of $\Gamma$ is given by sending a finite dimensional 
representation $\pi:\Gamma\to \End(\H)$ to its character 
$\chi_\pi(g)=\trace(\pi(g))$ (the \emph{non}-normalized trace on $\End(\H)$).  
If a class in $\KK^\Gamma(\C,\C)=\repgamma$
is represented by a $\Gamma$-invariant Fredholm operator $F: \H^{\ev}\to \H^{\odd}$,
then the corresponding virtual character in $\repgamma$ is given by the difference 
function 
$\chi=\chi_+-\chi_-$, where $\chi_+$ and $\chi_-$ denote the normalized traces 
of the $\Gamma$-representations  on $\H_+=\ker(F)$ and $\H_-=\coker(F)$, respectively. 
Since the value at a point $g\in \Gamma$ only depends on the action of $g$
on these spaces, it follows that in order to compute it we may always restrict
our attention to the closed subgroup $\Gamma_g\subseteq \O(V)$ generated by 
$\rho(g)$.  

Recall also that the identification $\KK^\Gamma(\C,\C)\cong \repgamma$ is multiplicative in the 
sense that it sends the Kasparov product $\hot_{\C}$ on $\KK^\Gamma(\C,\C)$ to the pointwise 
product of characters in $\repgamma$.

(b) We may always assume that $A\in \O(V)$. Indeed, if $A=O|A|$ is the 
polar decomposition of $A$ with $O=A|A|^{-1}$, then the homotopy 
$t\mapsto O(t\Id+(1-t)|A|)$ between $A$ and $O$ induces a $\Gamma$-invariant 
homotopy between the vector fields $A\rX$ and $O\rX$, and the result follows from 
Lemma \ref{lem-homotopy}.

(c) In case where $A= \Id$ is the identity, we obtain the class 
$[\rX]\in \K_0^\Gamma(C_\tau(V))$. It is the ``canonical'' generator of 
$\K_0^\Gamma(C_\tau(V))$ as described by Kasparov in \cite[\S 5]{Kas-operator}
and it follows from Kasparov's Bott-periodicity theorem \cite[Theorem 7 of \S 5]{Kas-operator}
that $[\rX]\otimes_{C_\tau(V)}[D]=1\in \repgamma$ (in the language of \cite{Kas-operator},
the class $[\rX]$ is denoted $\beta_V$ and $[D]$ is denoted $\alpha_V$).

Note that in case of the trivial group $\Gamma=\{e\}$  the above theorem
reduces to an index  computation given by L\"uck and Rosenberg in \cite{LR2}.
\end{remark}

The quantity we will be interested in for our Lefschetz theorem is the component of the trivial
representation in $\chi_{(\rho, A)}$: this is obtained by averaging the character over $\Gamma$; 
thus we derive the formula
\begin{corollary}
\label{fixedindex}
Suppose that $F:\H^{\ev}\to \H^{\odd}$  is a $\Gamma$-equivariant Fredholm operator 
representing  the Kasparov  product $[A\rX]\otimes_{C_\tau(V)}[D]\in \repgamma$. Then 
\begin{equation}
\mathrm{dim}_\C\, ( \mathrm{ker}^{\Gamma} F) -\mathrm{dim}_\C \,( \mathrm{coker}^{\Gamma}F)= 
\int_{\Gamma} \chi_{(\rho, A)}(g)\,dg
\end{equation}
(normalized Haar measure) where $V^{\Gamma}$ denotes the $\Gamma$-fixed points of a 
$\Gamma$-module $V$.
\end{corollary}

\begin{example}\label{ex-dim1}
In this example we want to compute the class 
$[A\rX]\hot_{C_\tau(\R)}[D_{\R}]\in \repgamma$ in the special case 
where $V=\R$ is one-dimensional, following the lines of Kasparov's \cite[Example 3 on p.760]{Kas0}. This gives the 
key calculation for the proof of Theorem \ref{charactertheorem}.
By part (b) and (c) of
the above remark we may assume that $A$ is multiplication by $-1$.
Also, by Lemma \ref{lem-support} we may restrict everything to the 
intervall $(-\pi,\pi)$.  If we identify $L^2(\Lambda^0_\C(-\pi,\pi))$
and $L^2(\Lambda^1_\C(-\pi,\pi))$ with $L^2(-\pi,\pi)$ in the canonical way,
we can realize the class $D=D_{(-\pi,\pi)}$ by the matrix $D=\left(\begin{matrix} 0 & \frac{d}{dx}\\
-\frac{d}{dx} & 0\end{matrix}\right)$. On the basis 
$\{e_n: n\in \Z\}$ 
with $e_n(x)=e^{inx}$ the operator $\frac{d}{dx}$ acts by $e_n\mapsto in e_n$, thus 
we obtain a bounded version $\tilde{d}: L^2(-\pi,\pi)\to L^2(-\pi,\pi)$ of the operator $d=\frac{d}{dx}$ by defining
$$\tilde{d} e_n= i\sign(n)e_n\quad \text{with}\;\sign(n)=\left\{\begin{matrix} 0&\text{if $n=0$}\\
\frac{n}{|n|}& \text{if $n\neq 0$}\end{matrix}\right\}.$$
The vector field $x\mapsto -x$ on $(-\pi,\pi)$  is homotopic to
$x\mapsto -\sin(\frac{x}{2})$. Thus, using Lemma \ref{lem-support} and the 
formula for the Kasparov product as given in (\ref{Kasp-product}) it follows that 
$[-\rX]\hot_{C_\tau(\R)}[D_\R]$ is given by the $\Gamma$-equivariant index
of the operator
$$T:=-\sin(\frac{x}{2})+\cos(\frac{x}{2})\tilde{d}:
L^2(-\pi,\pi)\to L^2(-\pi,\pi).$$
To compute it we first compute the index of the operator
$$S:=2i e^{i\frac{x}{2}}T= (1-e^{ix})+i(e^{ix}+1)\tilde{d}$$ which in 
terms of the orthonormal basis $\{e_n: n\in \Z\}$ 
is given by
$$S e_n=\left\{\begin{matrix} -2e_{n+1} & \text{if $n>0$}\\
e_0-e_1& \text{if $n=0$}\\
2 e_n &\text{if $n<0$}\end{matrix}\right\}.$$
It follows then from a short computation that $\ker S=\{0\}$ and $\coker S=< e_0+e_1>$.
Going back to the original operator $T$ we get
$\ker T=\{0\}$ and $\coker T$ is generated by
$\frac{1}{2}e^{-i\frac{x}{2}}(1+e^{ix})=\cos(\frac{x}{2})$.

If we write write $\O(\R)=\{1,-1\}$, 
then the corresponding action of $-1$ on 
$L^2(-\pi,\pi)\cong L^2(\Lambda^1_\C(-\pi,\pi))$ is given by $\xi\mapsto \big(x\mapsto -\xi(-x)\big)$.
Thus, on the generator $\xi(x)=\cos(\frac{x}{2})$ of $\coker T$ it acts by
multiplication with $-1$.
It follows that $[-\rX]\hot_{C_\tau(\R)}[D_{\R}]\in \repgamma$ is represented  by
the virtual character $\chi$ given by
$$\chi(g)=\left\{\begin{matrix} -1 & \text{if $\rho(g)=1$}\;\;\\
1 & \;\;\text{if $\rho(g)=-1$}\end{matrix}\right\}.$$
\end{example}

The following lemma will allow to reduce the proof of Theorem \ref{charactertheorem} to the case  
of the above example.

\begin{lemma}\label{lem-product}
For $i=1,2$ let $V_i$ be an Euclidean vector space with representation
$\rho_i:\Gamma\to O(V_i)$ and let $A_i\in \GL(V_i)$ commute with $\rho_i$.
Let $V=V_1\oplus V_2$, $\rho=\rho_1\oplus \rho_2$ and $A=A_1\oplus A_2$. Then 
$$[A\rX]\hot_{C_\tau(V)}[D_V]=
([A_1\rX_1]\hot_{C_\tau(V_1)}[D_{V_1}])\cdot([A_2\rX_2]\hot_{C_\tau(V_2)}[D_{V_2}]) \in \repgamma.$$
\end{lemma}
\begin{proof}
It is not difficult to check that under the canonical  isomorphism $C_\tau(V)\cong 
C_\tau(V_1)\hot C_\tau(V_2)$ we get 
$[A\rX]=[A_1\rX_1]\hot_{\C}[A_2\rX_2]$ in $\KK^\Gamma(C_{\tau}(V),\C)$ 
(compare with  the formula for $\beta_V$ in \cite[p. 546]{Kas-operator})
and it is shown in \cite[p. 547]{Kas-operator} that
$[D_V]=[D_{V_1}]\hot_{\C}[D_{V_2}]$ in $\KK^\Gamma(\C, C_\tau(V))$.
The result then follows from the associativity of the Kasparov product.
\end{proof}

\begin{proof}[Proof of theorem \ref{charactertheorem}]
Let $g\in \Gamma$ be fixed. As observed in Remark \ref{rem-character} 
we may assume that $\Gamma=\Gamma_g$ is the closed subgroup of $\O(V)$ 
generated by $\rho(g)$ (which we then identify with $g$). We also observed that
we may assume without loss of generality that $A\in \O(V)$.
Let $F\subseteq V$ be the set of $g$-fixed-points in $V$ and let $N=F^{\perp}$.
Then $F$ and $N$ are both, $g$- and $A$-invariant, and therefore 
the result will follow from the above lemma if we can show that
\begin{equation}\label{eq-F}
\big([A_F\rX]\hot_{C_\tau(F)}[D_F]\big)(g)=\sign\det(A_F),
\end{equation}
where $A_F$ denotes the restriction of $A$ to $F$, and 
\begin{equation}\label{eq-N}
\big([A_N\rX]\hot_{C_\tau(F)}[D_N]\big)(g)=1.
\end{equation}
Since $\Gamma_g$ acts trivially on $F$, we  may choose an orthonormal basis 
$\{v_1, \ldots, v_l\}$ of $F$  and, up to homotopy, we may assume
that  $A_F$  is given with respect to this basis by 
$\left(\begin{matrix} \pm 1& 0\\ 0& I_{l-1}\end{matrix}\right)$. If the upper left entry is $1$
we have $[A_F\rX]=[\rX]$ and the result follows from Kasparov's Bott-periodicity theorem 
(see Remark \ref{rem-character} (c)). If  the upper left entry is $-1$, we apply the above 
lemma  to the decomposition $F=< v_1>\oplus < v_2,\ldots, v_l>$.
Since $A_F$ restricts to the identity on $< v_2,\ldots, v_l>$, this summand provides the factor
$+1$ to the character at $g$ and since $g$ acts trivially on $< v_1>$ it follows from 
Example \ref{ex-dim1} that the first summand provides the factor $-1$ to the character at $g$.
This verifies (\ref{eq-F}).

To verify (\ref{eq-N}) we first consider the $-1$ eigenspace $V_{-1}$ for the action of $g$ on $N$.
This is clearly $\Gamma_g$- and $A$-invariant, and we may consider the decomposition 
$N=V_{-1}\oplus V_{-1}^{\perp}$ of $N$ as in the lemma. If $B$ denotes the restriction 
of $A$ to $V_{-1}$ we may again assume, up to $\Gamma_g$-invariant homotopy, that
$B= \left(\begin{matrix} \pm 1& 0\\ 0& I_{k-1}\end{matrix}\right)$ with respect to
 a suitable orthonormal base
$\{w_1,\ldots, w_k\}$ of $V_{-1}$. Decomposing 
$$V_{-1}=< w_1>\oplus < w_2,\ldots, w_k>$$
the second summand provides the factor $1$ by Bott-periodicity and the summand $< w_1>$ 
provides also the factor $1$ by Example \ref{ex-dim1}, since $g$ acts via the flip on $\R w_1$.

We therefore may assume without loss of generality that the action of $g$ on $N$ does not 
have eigenvalues $1$ or $-1$.  If $N\neq \{0\}$ let $\lambda_t=\cos(t)+i\sin(t)$ be a 
complex eigenvalue for the action of $g$ on the complexification $N_{\C}=N+iN$ of $N$ and
let $V_t\subseteq N_{\C}$ denote the corresponding eigenspace. Again, since $A$ commutes with 
$g$, it follows that $V_t$ is $A$-invariant. Since $A$ is orthogonal (and hence it acts 
unitarily on the complex vector space $V_t$) there exists a
non-zero $A$-eigenvector $u=u_1+iu_2\in V_t$ for some eigenvalue $\lambda_s=\cos(s)+i\sin(s)$,
 $s\in [0,2\pi)$.
It follows then from basic linear algebra that if we choose $u$ to be a unit vector in $N_\C$,
then $\sqrt{2}u_1, \sqrt{2}u_2$ are orthogonal unit vectors in $N$ and then $g$ and $A$ act on 
the invariant subspace $< u_1, u_2>\subseteq N$ via rotation by the angles $t$ and $s$, respectively.
But then we can $\Gamma_g$-equivariantly  homotop the restriction of $A$ to $< u_1, u_2>$ 
to the identity, which shows that the direct summand $< u_1, u_2>$ provides 
the factor $1$ to the character at $g$. Equation (\ref{eq-N}) now follows from a straightforward induction 
argument.
\end{proof}

In the remaining part of this section we want to discuss briefly
how Theorem \ref{charactertheorem} can also be obtained by
appealing to Atiyah and Singer
\cite{AS1}.  For ease of notation let $V=\R^n$ with standard inner poduct.
The cycle
$(L^2(\Lambda_{\C}^*\R^n),D+A\mathrm{X})$ is an unbounded
representative for the Kasparov product of the classes $[\rX]\otimes_{C_\tau(\R^n)}[D]$
 (see \citemore{Kas2}{Lemma 4}
and also \cite{Baaj-Julg} and \cite{Kuc} for the realization 
of $\KK$-classes by unbounded operators), it therefore is a
$\Gamma$-equivariant Fredholm operator on
$L^2(\Lambda^*_\C(\R^n))$ and has a $\Gamma$-equivariant index
$\indagamma (D+\mathrm{X}) \in \repgamma$ such that
$$
\indagamma(D+\mathrm{X}) = [\mathrm{X}]\hot_{\ctau (\Rn)} [D] \in \repgamma.
$$
We now eliminate Clifford algebras
from the picture, using the tangent bundle instead, using the well known 
$\KK^\Gamma$-equivalence between
$\ctau(\R^n)$ and $C_0(T\R^n)$ (a consequence of Kasparov's Bott-periodicity -- see 
\citemore{Kas-operator}{\S 5, Theorem 8},
\citemore{Blackadar}{24.5}). Under this equivalence
$[D]$ becomes the class $[\Dslash]$ of the \emph{Dolbeault
operator} on $T\Rn \cong \C^{n}$, and $[\rX]$ becomes in the
notation of Atiyah-Singer the Bott generator,
$j_0!(1)\in \KK^{\Gamma}(\C, C_0(T\Rn)) =
\K^{0}_{\Gamma}(T\Rn)$, where $j_0\colon \{0\}\to \Rn$ is the
inclusion of the origin of $\Rn$. Atiyah
and Singer say the index map takes the class $j_0!(1)$ to $1$.
On the other hand, the class
$[A\mathrm{X}]\in \K_0^{\Gamma}(\ctau(\Rn))$  corresponds
to $A_*(j_0!(1))\in \K^0_{\Gamma}(T\Rn)$. Therefore, following
\cite{AS1}, computing
$[A\mathrm{X}]\hot_{\ctau(\Rn)}[D]=\indagamma(D+A\mathrm{X})$ is
equivalent to computing the topological index 
 $\indtgamma\bigl(A_*(j_0!(1))\bigr)\in
\repgamma$ as introduced in \cite{AS1}.

As before, in order to compute the character at $g$
 it suffices  to assume that $\Gamma = \Gamma_g$ is the subgroup of $\on$
 generated by $\rho(g)$. 
For convenience, let $\beta_W \in \kgamma (TW)$ be the Bott
generator, whenever $W$ is a $\Gamma$-invariant linear subspace
of $\Rn$. As is well-known, $\kgamma (TW)$ is a rank-one
$\repgamma$-module with generator $\beta_W$. Equivariant Bott
periodicity $\indtgamma\colon \kgamma (TW) \to \repgamma$
commutes with the module action, and $\indtgamma (\beta_W) = 1
\in \repgamma$ by Atiyah-Singer \cite{AS1}. To be explicit, let
$\sigma \colon \pi^{*}E \to \pi^{*}E$ be an odd endomorphism of
$\Z/2$-graded bundles, with $\sigma$ an isomorphism outside of
a compact subset of $T\Rn$, and so representing a class $a\in
\kgamma (T\Rn)$. Suppose $b \in \repgamma$ is represented by a
finite-dimensional $\Gamma$-vector space $V$. Then $a\cdot b$
is represented by $\sigma\otimes \mathrm{id}_V\colon
\pi^{*}(E\otimes V) \to \pi^{*}(E\otimes V)$.

The cycle for $\kgamma (T\Rn)$ representing $A_{*}\beta_{\Rn}$
is given by the trivial $\Z/2$-graded bundle $T\Rn \times \ext
(\Rn)$ together with the odd endomorphism $\sigma\colon T\Rn
\times \ext (\Rn) \to T\Rn \times \ext (\Rn)$ determined by the
map $T\Rn \to \C^n$, $(x,\xi) \mapsto Ax + i\xi$ (using
Clifford multiplication.) Note that as $Ax +i\xi$ vanishes only
at the origin of $T\Rn$, the endomorphism $\sigma$ is an
isomorphism outside of a compact set.

Let $F $ denote the fixed subspace of $g$ and $N = F^{\perp}$. We have a well-known
isomorphism
\begin{equation}
\label{tensorproduct}
 \ext (\Rn) \cong \ext (F)\hot \ext (N)\end{equation}
of graded vector spaces, and there is a corresponding isomorphism of (trivial) bundles.
Note that $F$ and $N$ are also $A$-invariant.

If we restrict $\sigma\colon T\Rn \times \ext (\Rn) \to T\Rn \times \ext (\Rn)$ to $TF$, then
under the identification \eqref{tensorproduct},
the endomorphism $\sigma$, when restricted to $TF$, becomes
the endomorphism $\sigma\hot \mathrm{id}_{N} \colon TF\times \ext (F)\hot \ext (N) \to TF \times \ext (F)\hot \ext (N).$
Thus we have the following.

\begin{lem}
If $i_F \colon F \to \Rn$ is the $\Gamma$-equivariant inclusion, then
\begin{equation}
i^{*}_{F}(A_{*}(\beta_{\Rn})) = \mathrm{sign}\, \mathrm{det} (A_{|_{F}}) \,  \beta_F  \cdot [\ext (N)]\in  \kgamma (TF),
\end{equation}
 where $\ext (N) \in \repgamma$ is given by
$\sum_{i=0}^{\mathrm{dim}(N)} (-1)^{i} [\Lambda^{i}_\C N]$, an alternating sum of finite-dimensional
$\Gamma$-spaces.
\end{lem}
\begin{proof}
Given the preceding discussion, it is clear that
\begin{equation}
i^{*}_{F}(A_{*}(\beta_{\Rn})) = (A_{|_{F}})_{*}(  \beta_F)  \cdot [\ext (N)]\in  \kgamma (TF).
\end{equation}

Since $\Gamma$ acts trivially on $F$ (and $TF$), if the
restriction of $A$ to $F$ has positive determinant, it is
$\Gamma$-equivariantly
homotopic  to
the identity. If the determinant is negative, it is similarly
homotopic to a reflection $Q$, and it is standard that
 $Q_{*}(\beta_F) = -\beta_F$ in non-equivariant $\K$-theory, but then in
this case also, because the $\Gamma$-action on $TF$ is trivial.
\end{proof}

Following a pattern of argumentation in Atiyah-Segal
\cite{AS2}, since $\rho(g) $ the generator of $\Gamma$ has no fixed
points in $N$, the class $[\ext (N)]$ is a unit in the
localization $\repgamma_g$ of the ring $\repgamma$ at the prime
ideal determined by $g$ (see \cite[Lemma 2.7]{AS2})--indeed, this prime ideal consists of
all characters which vanish at $g$, while the character
corresponding to $[\ext (N)]$ is 
$$g \mapsto
\sum_{i=0}^{\mathrm{dim}(N)} (-1)^{i} \mathrm{trace}(g\colon
\Lambda^{i}_\C N \to  \Lambda^{i}_\C N) =
 \mathrm{det}(1-g_{|_{N}}) \not= 0.$$
  For this reason and the above calculation, we see that
$i_F^{*}\colon \kgamma (T\Rn) \to \kgamma (TF)$ is an
isomorphism after localizing at $g$ (c.f. \cite[Proposition 2.8]{AS2}). 
Since $i_F^{*}(\beta_{\Rn}) =
\beta_{F}\cdot [\ext (N)]$ by the same argumentation with $A$
set equal to the identity, we get that $i_F^{*}(
A_{*}\beta_{\Rn}) = \mathrm{sign}\, \mathrm{det} (A_{|_{F}})
\,i_F^{*}(\beta_{\Rn})$ and hence since $i_F^{*}$ is an
isomorphism after localization at $g$, that $A_{*}(\beta_{\Rn})
= \mathrm{sign}\, \mathrm{det} (A_{|_{F}}) \beta_{\Rn}$ after
localization at $g$. Therefore, taking $\indtgamma$ of both
sides and using that $\indtgamma (\beta_{\Rn}) = 1 \in
\repgamma$, gives that $\indtgamma (A_{*}\beta_{\Rn}) =
\mathrm{sign}\, \mathrm{det} (A_{|_{F}})\, 1 \in \repgamma_g$.
Evaluation of characters at $g$ passes  to the
localization, and is compatible with evaluation before
localization, whence evaluating the above expression at $g$
gives that $\indtgamma (A_*\beta_{\Rn})(g) = \mathrm{sign}\,
\mathrm{det} (A_{|_{F}})$ as required. This gives the second proof of Theorem
\ref{charactertheorem}.

\section{The Lefschetz theorem}\label{lefschetz}

Let the countable group $G$ act isometrically, properly and
co-compactly on the Riemannian manifold $X$ (it follows that $X$ is
complete.) Let $\phi\colon X \to X$ be a smooth map. We are going to
formulate and prove a Lefschetz fixed-point formula in this context
using the discussion in Section \ref{fundamental} on Poincar\'e
duality between $\Lambda=C_0(X)\rtimes G$ and
$\Lambdah=C_{\tau}(X)\rtimes G$.
To get an endomorphism of the algebra
$\Lambda$ and to be adequate for the formulation of the Lefschetz
theorem, we need a couple of assumptions on the map $\phi$ on the
manifold $X$.

First we require a transversality of $\phi$. Suppose for the moment
that the $G$-action on $X$ is free. Then $G\backslash X$ is a
manifold, and since $\phi$ maps orbits to orbits, we obtain a smooth
map $\dot{\phi}\colon G\backslash X \to G\backslash X$. In this
case, we want to demand that $\dot{\phi}$ is in \emph{general
position}: that is, that its graph is transverse to the diagonal in
$G\backslash X \times G\backslash X$.

By definition of the smooth structure on $G\backslash X$, this means
the following: \emph{If $x \in X$, $g\in G$ such that $\phi(gx) =
x$, then the map
\begin{equation}\label{eq-gp}
\mathrm{Id} -d(\phi \circ g) (x)\colon T_xX \to
T_xX
\end{equation}is non-singular.}

If the $G$-action is not free, $G\backslash X$ is not a manifold.
But the reformulation of the condition that $\dot{\phi}$ be in general
position given above still makes sense. We thus impose the
following:

\begin{assumption}\label{A1}For every $g \in G$, the smooth map $\phi \circ g\colon X \to
X$ is in general position, i.e., (\ref{eq-gp})  holds for all $x\in X$ with $\phi(gx)=x$.
\end{assumption}

Next we require the following compatibility of the map $\phi$ and
the $G$-action on $X$:

\begin{assumption}\label{A2} There is an automorphism $\zeta\colon G \to G$
such that
\begin{equation}\label{compatibilitycondition}
\phi \bigl( \zeta (g)x\bigr) = g\bigl( \phi (x)\bigr), \; \text{for
all $x\in X$}.
\end{equation}
\end{assumption}

The assumption ensures that the maps $f\mapsto f\circ \phi$ and $g
\mapsto \zeta (g)$ constitute a covariant pair for the action of $G$
on $C_0(X)$. We obtain an automorphism
\begin{equation}\label{alpha}\alpha\colon C_0(X)\cross G \to
 C_0(X)\cross G.
 \end{equation}

The abstract Lefschetz theorem (see \cite{Emerson}) asserts
that the Leftschetz number equals an index-theoretic pairing,
i.e.,
\begin{equation}\label{abtractLef}
\mathrm{Lef}([\alpha])=\langle\widehat{[\alpha]},\Delta\rangle,
\end{equation}
where $\mathrm{Lef}([\alpha])$ is the Lefschetz number of
$\alpha$ defined by
$$\mathrm{Lef}([\alpha]) 
\defeq \mathrm{tr}_s (\alpha_*\colon \K_*(C_0(X)\cross G)_{\QQ} \to \K_*(C_0(X)\cross
G)_{\QQ})
$$
(which only depends on $[\alpha] \in \KK(C_0(X)\cross G,
C_0(X)\cross G)$) and $\widehat{[\alpha]}$ denotes the
Poincar\'e dual of $[\alpha]$ which more exactly equals
 $\alpha_*({\sigma_{}}_*\Deltah)\in \KK(\C, \Lambda \hot \Lambdah)$, where, as before,
 $\sigma: \Lambdah\hot\Lambda\to \Lambda\hot\Lambdah$ denotes the flip map.

Therefore to prove Theorem \ref{intro:mainthm}, we want to
compute
 the pairing
$\alpha_*({\sigma_{}}_*\Deltah)\hot_{\Lambda\hot \Lambdah}\Delta,$
where $\alpha$ is as in \eqref{alpha}. By functoriality, this is the
same as ${\sigma_{}}_*\Deltah\hot_{\Lambda\hot \Lambdah}
\alpha^{*}(\Delta)$, which we will focus on instead.

We set
\begin{equation}
F_\epsilon \defeq \{(x,g) \in X\times G\mid \rho\bigl( \phi(gx),
x\bigr)< \epsilon\},
\end{equation}
where $\rho$ denotes the metric on $X$.
Give $F_\epsilon$ the structure of a $G$-space by restricting the
following action of $G$ on $X\times G$:
\begin{equation}
\label{gaction}
h(x,g) \defeq (hx, \zeta (h)gh^{-1}).
\end{equation}
Let $F = F_0$ in the above notation, so $F = \{(x,g)\mid \phi(gx) =
x\}$. Then $F_\epsilon$ is a neighbourhood of $F$ and $F_\epsilon
\to F$ as $\epsilon \to 0$. Note also that $G$ leaves $F_\epsilon$
(and likewise $F = F_0$) invariant, as if $\rho\bigl( \phi(gx), x\bigr)
< \epsilon$ then
$$ \rho\bigl( \phi (\zeta (h)gh^{-1}hx), hx\bigr) = \rho\bigl( h \phi (gx), hx\bigr) = \rho\bigl( \phi (gx),x\bigr) < \epsilon. $$
Let $V_\epsilon$ be the set of first coordinates of points in
$F_\epsilon$. Then $V_\epsilon$ is a $G$-set for $\epsilon \ge 0$.
Let $V\defeq V_0$.

\begin{lem}
\label{components}
The set $V$ is discrete. Furthermore, if $\delta > 0$, there exists
$\epsilon > 0$ such that every component of $V_\epsilon$ is
contained in a $\delta$-ball in $X$ with center in $V$.
\end{lem}
\begin{proof}
Suppose $(x_j)$ and $(g_j)$ are sequences in $X$ and $G$
respectively such that $\phi(g_jx_j) = x_j$, the $x_j$ are all
distinct, and $x_j \to x_0$ for some $x_0$. Let $h_j$ such that
$\zeta(h_j) = g_j^{-1}$. Then $\phi(x_j) =  h_jx_j$. Since $x_j \to
x_0$, $\phi(x_j) \to \phi(x_0)$, and hence $h_jx_j \to \phi(x_0)$.
But then
$$\rho\bigl(h_jx_0, \phi (x_0)\bigr) \le \rho( h_jx_0, h_jx_j) + \rho\bigl(h_jx_j, \phi (x_0)\bigr) \to 0.$$
But since the $G$-action is proper, there are only finitely many $h
\in G$ which map $x_0$ to any fixed, pre-compact neighbourhood of
$\phi(x_0)$. Hence $h_j = h$ for some $h$ and almost all $j$. We may
assume $h_j = h$ for all $j$, which gives that $g_j = g$ for all $j$
and then $\phi \circ g$ has an accumulation point amongst its fixed
points, which contradicts Assumption \ref{A1}. This argument proves that
$V$ is discrete.

For the second statement, observe that $\{V_{\frac{1}{n}}\}$ is a
nested sequence whose intersection is $V$. Using the $G$-compactness
of $X$, we see that, given $\delta>0$, there exists $n\in\N$ such that
$V_{\frac{1}{n}}$ is contained in the
$\delta$-neighbourhood of $V$. By the first
statement, the second statement now follows.
\end{proof}

Since $V$ is discrete and the $G$-action on $X$ is co-compact, $V$
splits into finitely many $G$-orbits. Observe that the set of such
orbits has an obvious correspondence with the set
\begin{equation}\label{Fix}\operatorname{Fix}(\dot{\phi}):=\{p\in
X\,|\,\dot{\phi}(\dot{p})=\dot{p}\},
\end{equation}
where $\dot{\phi}$ is the induced map $G\backslash X\to G\backslash
X$ and $\dot{p}$ denotes an orbit of $p$. Let us denote each
$G$-orbit in $V$ corresponding to each point $p\in
\operatorname{Fix}(\dot{\phi})$ by $V_p$.

The $G$-set $F$ admits a similar decomposition, $F = \sqcup F_p$,
where $F_p = \{(x,h) \in F \mid x \in V_p\}$. For each $V_p$ fix an
element $g_p\in G$ such that $\phi(g_pp) = p$. Let $L_p
\defeq g_pK_p$ be the coset of $K_p\defeq \mathrm{Stab}_G(p)$. Then
one can see that $L_p=\{g\in G\, |\, \phi(gp)=p\}$.

From this, we get the following. Consider a point $gp \in V_p$. Then
there exists $h \in G$ such that $\phi(hgp) = gp$, and hence
$\phi\bigl( \zeta(g)^{-1}hgp\bigr) = p = \phi(g_pp)$, so that $g_pp
= \zeta (g)^{-1}hgp$ and $g_p^{-1}\zeta (g)^{-1}hg\in K_p$.
\emph{Hence $h$ lies in the twisted conjugate $\zeta (g)L_p g^{-1}$
of $L_p$.} The converse of this statement is also true.

Hence we can write
$$V_p = \{gp \mid gK_p \in G/K_p\}, \;\;\; F_p = \{ (gp, h) \mid gK_p \in G/ K_p , \; h \in \zeta (g)L_pg^{-1}\}.$$

Similarly, we get decompositions of $V_{\epsilon}$ and
$F_{\epsilon}$. By Lemma \ref{components}, we may choose
$\delta>0$ small enough so that all the $\delta$-balls centered
at the points of $V$ are disjoint and therefore there exists
$\epsilon>0$ such that
\begin{equation}\label{Vepsilon} V_{\epsilon} = \sqcup_{p\in
\operatorname{Fix}(\dot{\phi}), gK_p\in G/K_p} \;
V_{\epsilon,gp},
\end{equation} where $V_{\epsilon,gp}$ is the part of
$V_{\epsilon}$ which is contained in the $\delta$-ball centered at the
point $gp\in V$. Similarly,
\begin{equation}\label{Fepsilon}
F_{\epsilon} = \sqcup_{p\in\operatorname{Fix}(\dot{\phi}), gK_p \in G/K_p} \; F_{\epsilon,gp},
\end{equation}
where
$$F_{\epsilon,gp} = \{ ( x, h) \in X\times G \mid x\in V_{\epsilon,gp},\; h \in \zeta(g)L_pg^{-1}\}.$$

In what follows, we shall decribe the pairing
${\sigma_{}}_*\Deltah\hot_{\Lambda\hot
\Lambdah}\alpha^*\Delta$ as a direct sum of Kasparov products 
which live on the Hilbert spaces $L^2(\Lambda^*_{\CC}(V_{\eps,p}))^{\Gamma_{p,g}}$,
where $\Gamma_{p,g}\subseteq K_p$ denotes the stabilizer of $g\in L_p$ under 
the conjugation action $g\mapsto \zeta(h)gh^{-1}$. These summands can then be computed
via the results of the previous section.
We start with a careful description of the Hilbert
space (recall Definitions \ref{classdelta} and
\ref{classdeltah}). It is the tensor product of the right
Hilbert $\Lambda\hot \Lambdah$-module $\mathcal{E}$ described
prior to Definition \ref{classdeltah}, and the Hilbert space
$\ltwoformsg$ occurring in connection with the fundamental
class $\Delta$, twisted by the automorphism $\alpha$ induced
from $\phi$ and $\zeta$ (see \eqref{alpha}). Here and throughout we write
$\Lambda$ for $C_0(X)\rtimes G$ annd $\Lambdah$ for $C_\tau(X)\rtimes G$.

After twisting $\Delta$ by $\alpha$, we obtain the Hilbert space
$\ltwoformsg$ equipped with a twisted representation of $\Lambda\hot
\Lambdah$ whose explicit form we state for the record (compare with the
untwisted version in \eqref{representationCtau} and
\eqref{representationCzero}): the algebra $\Lambda=C_0(X)\cross G$ acts via the covariant pair
\begin{equation}\label{therep}
f(\xi\hot e_g) = g^{-1}(f\circ \phi)\xi\hot e_g, \quad\quad
h\cdot(\xi\hot e_g)= \xi\hot e_{\zeta (h)g},
\end{equation}
for $ f\in C_0(X),h\in G$ and $\xi\hot e_g\in \ltwoformsg$.
The algebra  $\Lambdah=C_\tau(X)\cross G$ acts 
by the covariant
pair
\begin{equation}\label{therep1}
\varphi (\xi\hot e_g) = \varphi \xi \hot e_g, \quad\quad\quad h\cdot(\xi\hot
e_g) = h(\xi)\hot e_{gh^{-1}},
\end{equation}
for $\varphi\in \ctau(X), h\in G$ and $\xi\hot e_g\in \ltwoformsg$.
Recall that $\mathcal{E}$ is the completion of
$C_c(X){\F_U}\hot \C G$ with respect to a certain inner
product, where $U=\{(x,y)| \rho(x,y)<\epsilon\}$ for some
$\epsilon>0$ from Remark \ref{remark_Kasparov} (iv). We may choose (and fix)
 $\epsilon$ such that 
(\ref{Vepsilon}) is satisfied for a suitable $\delta>0$. 

Notice that there is a well defined inclusion of the algebraic tensor product
$C_c(X)\F_U\odot \C G$ into $\Lambda\hot\Lambdah$
given by sending the elementary tensor $\alpha\hot [h]$ to the element
$\alpha([h]\hot [h])\in 
\big(C_0(X)\hot C_\tau(X)\big)\rtimes (G\times G)\big)\cong \Lambda\hot\Lambdah$.
Using this inclusion, we obtain a natural pairing
$$M : \big(C_c(X)\F_U\odot \C G\big)\times \big(\ltwoformsg\big)\to \ltwoformsg $$
given by  applying the action of $\Lambda\hot\Lambdah$ as described in (\ref{therep1})
to the image of $C_c(X)\F_U\odot \C G$ under the above described inclusion.

Recall that $F_{\eps}=\{(x,g)\in X\times G\; |\; \rho(\phi(gx),x)<\eps\}$. We denote by 
$L^2(\Lambda^*_{\CC}F_{\eps})$ the set of all $\xi\in \ltwoformsg$ which \emph{live}
on $F_\eps$ in the obvious sense (by viewing  the elements of $\ltwoformsg$ as sections
on $X\times G$). We then get

\begin{lemma}\label{lem-action}
Let $\alpha\hot[h]\in C_c(X)\F_U\odot \C G$ act on $\ltwoformsg$ as described above.
Then $(\alpha\hot[h])\cdot \ltwoformsg\subseteq L^2_c(\Lambda^*_{\CC}F_{\eps})$, where 
$L^2_c(\Lambda^*_{\CC}F_{\eps})$ denotes the set of $L^2$-sections on $F_{\eps}$ which
vanish outside some compact subset of $F_{\eps}$.
\end{lemma}
\begin{proof}
If we regard the elements of $\ltwoformsg$ as sections on $X\times G$ in the canonical way,
it follows from  (\ref{therep1}) that the action of $\alpha\hot [h]$
on such section $\mu\in \ltwoformsg$ is given by the formula
\begin{equation}\label{eq-M}
\big((\alpha\hot[h])\cdot \mu\big)(x,g)=\alpha(\phi(gx),x)d_h^x(\mu(h^{-1}x, \zeta(h^{-1})gh))),
\end{equation}
where $d_h^x:\Cl(T_{h^{-1}x}X)\to \Cl(T_xX)$ is the isomorphism underlying the 
action of $G$ on $\ctau(X)$.
Thus the result follows directly  from the fact that $\alpha$ is compactly supported in 
$U_{\eps}=\{(x,y) : \rho(x,y)<\eps\}$.
\end{proof}

In what follows let $P\subseteq X$ be a fixed set of representatives for
 $\Fix(\dot\phi)=G\backslash V$, with $V$ as in the discussion in the beginning 
 of this section. Let 
 $$S:= \{(p, gK_p): p\in P, gK_p\in G/K_p\},$$
  where 
 $K_p$ denotes the stabilizer of $p$ in $G$. 
  Recall from (\ref{Fepsilon}) that  for $\eps>0$ small enough, the set $F_{\eps}$ decomposes into a disjoint union
$$F_{\eps}= \sqcup_{(p, gK_p)\in S} \; F_{\epsilon,gp},$$
with 
$$F_{\epsilon,gp} = \{ ( x, h) \in X\times G \mid x\in V_{\epsilon,gp},\; h \in \zeta(g)L_pg^{-1}\}.$$
 Note that $G$ acts unitarily on $L^2(\Lambda^*_{\CC}F_{\eps})$ via 
 \begin{equation}\label{actionFeps}
 (s\mu)(x,g)=d_s^x\big(\mu(s^{-1}x, \zeta(s^{-1})gs)\big),
 \end{equation}
 for $s\in G$, $\mu\in L_c^2(\Lambda^*_\C F_{\eps})$,
 where, by abuse of notation, $d_s^x:\Lambda_\C^*(T_{s^{-1}x}X)\to \Lambda_\C^*(T_{x}X)$ 
 is the isometry induced by the differential $d_s^x:T_{s^{-1}x}X\to T_xX$.
 This action restricts to well defined actions of $K_p=\Stab_G(p)$ 
on $F_{\eps,p}$ for all $p\in P$. As usual, we let $L^2(\Lambda^*_\C F_{\eps,p})^{K_p}$
denote the $K_p$-invariant elements in $L^2(\Lambda^*_\C F_{\eps,p})$.

In what follows, we equip $L^2_c(\Lambda^*_{\CC}F_{\eps})$ with a new inner product 
given by
$$\lk \mu, \nu\rk=\sum_{s\in G} \lk s(\mu),\nu\rk_{L^2(\Lambda^*_{\CC}F_{\eps})}$$
with action of $G$ on $L^2_c(\Lambda^*_{\CC}F_{\eps})$ as explained above.
Note that this inner product makes sense, since $G$ acts properly on $F_{\eps}$ 
and $\mu$ and $\nu$ are compactly supported. We denote by $\H_\eps$ the Hausdorff
completion 
of $L^2_c(\Lambda^*_{\CC}F_{\eps})$ with respect to this inner product.

\begin{lemma}\label{lem-hilbert} 
Consider the composition $\Phi=\Psi\circ M$ of maps
$$
\begin{CD}
\E\hot_{\Lambda\hot\Lambdah} \big(\ltwoformsg\big)@>M>> \H_\eps
@>\Psi>> \oplus_{p\in P} L^2(\Lambda^*_\C F_{\eps,p})^{K_p},
\end{CD}
$$ 
where $M$ is given on elementary tensors by the pairing of Lemma \ref{lem-action}
and where 
$$\Psi(\eta)=\oplus_{p\in P} \frac{1}{\sqrt{|K_p|}}\sum_{s\in G} s\eta|_{F_{\eps,p}},$$
for $\eta\in L^2_c(\Lambda^*_{\CC}F_{\eps})$.
Then $\Phi$ is an isometric isomorphism of Hilbert spaces.
\end{lemma}
\begin{proof}
Using formulas (\ref{innerproduct}), (\ref{therep}), 
(\ref{therep1}) and (\ref{eq-M}),
we compute for all $\alpha_i\hot[h_i]\in C_c(X)\F_U\hot\C G$ and 
$\mu_i\in \ltwoformsg$, $i=1,2$:

\begin{align*}
&\langle (\alpha_1\hot [h_1])\hot_{\Lambda\hot\Lambdah}\mu_1,\,(\alpha_2\hot [h_2])\hot_{\Lambda\hot\Lambdah}\mu_2\rangle_{\mathcal{E}\hot_{\Lambda\hot \Lambdah}
\ltwoformsg}\\
&=\sum_{g\in G}\int_X
\langle \big(\langle \alpha_2\hot [h_2],\, \alpha_1\hot [h_1]\rangle_{\E} \cdot\mu_1\big)(x,g),\, 
\mu_2(x,g)\rangle\, dx\\
&=\sum_{g\in G}\sum_{s\in G}\int_X
\langle \big( h_2^{-1}\left(\alpha_2^* s(\alpha_1)\right)[h_2^{-1}sh_1]\hot [h_2^{-1}s]\big)\cdot\mu_1(x,g),\, \mu_2(x,g)\rangle\, dx\\
&=\sum_{g\in G}\sum_{s\in G}\int_X\lk h_2^{-1}\big(\alpha_2^*s(\alpha_1)\big)(\phi(gx), x) d_{h_2^{-1}sh_1}^x\big(\mu_1(h_1^{-1}s^{-1}h_2x, \zeta(s^{-1}h_2)gh_2^{-1}sh_1)\big),\\
&\quad\quad\quad\quad\quad\quad\quad\quad\quad\quad\quad\quad\quad\quad\quad\quad\quad\quad\quad\quad\quad\quad\quad\quad\quad\quad\quad\quad\quad\quad\quad\quad
\mu_2(x,g)\rk\, dx\\
&=\sum_{g\in G}\sum_{s\in G}\int_X
\lk d_{h_2^{-1}}^x\big(\alpha_2^*s(\alpha_1)(h_2\phi(gx), h_2x) \big)d_{h_2^{-1}sh_1}^x\big(\mu_1(h_1^{-1}s^{-1}h_2x, \zeta(s^{-1}h_2)gh_2^{-1}sh_1)\big),\\
&\quad\quad\quad\quad\quad\quad\quad\quad\quad\quad\quad\quad\quad\quad\quad\quad\quad\quad\quad\quad\quad\quad\quad\quad\quad\quad\quad\quad\quad\quad\quad\quad
 \mu_2(x,g)\rk\, dx
\end{align*}
Now, applying on both sides the unitary transformation $\nu\mapsto h_2\nu$  given by the formula
in (\ref{actionFeps}), the above term transforms into
\begin{align*}
&=\sum_{g\in G}\sum_{s\in G}\int_X
\lk \big(\alpha_2^*s(\alpha_1)\big)(\phi(gx), x) d_{sh_1}^x\big(\mu_1(h_1^{-1}s^{-1}x,\zeta(s^{-1})gsh_1)\big),\\
&\quad\quad\quad\quad\quad\quad\quad\quad\quad\quad\quad\quad\quad\quad\quad\quad\quad\quad\quad\quad\quad\quad\quad
d_{h_2}^x\big(\mu_2(h_2^{-1}x,\zeta(h_2^{-1})gh_2)\big)\rk\, dx\\
&=\sum_{g\in G}\sum_{s\in G}\int_X
\lk \big(d_s^x\big(\alpha_1)(s^{-1}\phi(gx), s^{-1}x) d_{sh_1}^x\big(\mu_1(h_1^{-1}s^{-1}x, \zeta(s^{-1})gsh_1)\big),\\
&\quad\quad\quad\quad\quad\quad\quad\quad\quad\quad\quad\quad\quad\quad\quad\quad\quad
\alpha_2(\phi(gx),x)d_{h_2}^x\big(\mu_2(h_2^{-1}x,\zeta(h_2^{-1})gh_2)\big)\rk\, dx\\
&=\sum_{s\in G}\sum_{g\in G}\int_X\langle s\big((\alpha_1\hot[h_1])\cdot\mu_1\big)(x,g),\, 
(\alpha_2\hot[h_2])\cdot\mu_2(x,g)\rangle\, dx\\
&=\sum_{s\in G} \lk s\big((\alpha_1\hot[h_1])\cdot\mu_1\big), (\alpha_2\hot[h_2])\cdot\mu_2\rk_{L^2(\Lambda^*_\C F_{\epsilon})}.
\end{align*}
This shows that $M$ extends to a well defined unitary homomorphism 
from $\E\hot_{\Lambda\hot\Lambdah} \big(\ltwoformsg\big)$ to $\H_\eps$, and it is not 
difficult to see that it has dense image. Thus the result will follow if we can show that
$\Psi: \H_\eps\to \oplus_{p\in P} L^2(\Lambda^*_\C F_{\eps,p})^{K_p}$ is also isometric 
(it clearly has dense image). After decomposing $F_{\eps}$ into the disjoint union
$\sqcup_{(p, gK_p)\in S} \; F_{\epsilon,gp}$, we may assume without loss of generality that
$P=\{p\}$ is a single point. Then $L^2_c(\Lambda^*_\C F_{\eps})$
can be written as the set of finite sums $\eta=\sum_{gK_p} \eta_g$ 
with $\eta_g$ supported on $F_{\eps, gp}$. Each such function is of the form $g\eta'$ for some
$\eta'\in L^2(\Lambda^*_\C F_{\eps,p})$. So assume now that 
$\eta, \nu\in  L^2(\Lambda^*_\C F_{\eps,p})$ and $g,h\in G$. Then
\begin{align*}
\lk g\eta, h\nu\rk_{\H_\eps}&=
\sum_{s\in G} \lk sg\eta, h\nu\rk_{L^2(\Lambda^*_\C F_{\eps})}
\stackrel{s\mapsto hsg^{-1}}{=} \sum_{s\in G} \lk s\eta, \nu\rk_{L^2(\Lambda^*_\C F_{\eps})}
=\sum_{s\in K_p}\lk s\eta,\nu\rk_{L^2(\Lambda^*_\C F_{\eps})}
\end{align*}
On the other hand, we have 
$$\Psi(g\eta)=\frac{1}{\sqrt{|K_p|}}\sum_{s\in G} sg\eta|_{F_{\eps,p}}=\frac{1}{\sqrt{|K_p|}}\sum_{s\in K_p} 
s\eta|_{F_{\eps,p}}$$
from which we get
\begin{align*}
\lk\Psi(g\eta), \Psi(h\nu)\rk_{L^2(\Lambda^*_C F_{\eps,p})}&=
\frac{1}{|K_p|}\sum_{s,t\in K_p}\lk s\eta, t\nu\rk_{L^2(\Lambda^*_C F_{\eps,p})}
=\sum_{s\in K_p}\lk s\eta, \nu\rk_{L^2(\Lambda^*_C F_{\eps,p})},
\end{align*}
which now proves that $\lk\Psi(\eta),\Psi(\nu)\rk_{L^2(\Lambda^*_C F_{\eps,p})}=
\lk \eta,\nu\rk_{\H_{\eps}}$ for all $\eta,\nu\in L^2_c(\Lambda^*_\C F_{\eps})$.
\end{proof}

Since $K_p$ acts on $F_{\epsilon, p}=V_{\epsilon,p}\times L_p$
by $h\cdot(x,g)=(hx,\zeta(h)gh^{-1})$ as defined in
\eqref{gaction}, and also since one can consider
$L^2(\Lambda^*_{\C}F_{\epsilon, p})$ as a direct sum of copies
of $L^{2}(\Lambda^*_{\C} V_{\epsilon,p})$, one summand for each
point in $L_p$, we have
$$L^2(\Lambda^*_{\C}F_{\epsilon,p})^{K_p}=\bigoplus_g
L^{2}(\Lambda^*_{\C}V_{\epsilon,p})^{\Gamma_{p,g}},$$ where $g$ runs through
a given set $\Sigma_p$ of representatives for the orbits in $L_p$ under the twisted
conjugation by $K_p$ and $\Gamma_{p,g}\subset K_p$ denotes the
stabilizer of $g$ under this action. Thus, combining this observation with
the above lemma we get

\begin{equation}\label{hilbertspace2}
\mathcal{E}\hot_{\Lambda\hot \Lambdah}
\big(\ltwoformsg\big)\cong\bigoplus_{\Sigma}
L^{2}(\Lambda^*_{\C}V_{\epsilon,p})^{\Gamma_{p,g}}=\bigoplus_{\Sigma} \H_{p,g}^{\Gamma_{p,g}},
\end{equation} with $\Sigma=\cup_{p\in P}\Sigma_p$
and $\H_{p,g}=L^{2}(\Lambda^*_{\C}V_{\epsilon,p})$.

We are now going to compute the operator. For this let
$g\in L_p$. Since $\rho(\phi(gx),x)<\eps$ for all $x\in V_{\eps,p}$ 
and $g\in L_p$, we have $(\phi(gx),x)\in U_\eps$ for all such $x$ and $g$.
Thus we have a well defined vector field $\theta_{p,g}:V_{\eps,p}\to 
TV_{\eps,p}$  given by $\theta_{p,g}(x)=\theta_{\eps}(\phi(gx),x)$ with 
$\theta_{\eps}(z,x)=\frac{\rho(z,x)}{\eps}d_x\rho(z,x)$ as 
in Remark \ref{remark_Kasparov} (iv). It determines a class 
$ \Theta_{g,p}\in \KK^{\Gamma_{g,p}}(\CC, C_\tau(V_{\eps,p}))$
as in the previous section. 
Indeed,
since $\|d_x\rho(z,x)\|=1$ for all $z,x\in X$ with $z\neq x$, 
it follows that $\theta_{p,g}(x)^2-1=\|\theta_{p,g}(x)\|^2-1\to 0$ if $x\to \infty$ in $V_{\eps,p}$,
and therefore the class $\Theta_{p,g}$ is given directly via Clifford multiplication 
of $\theta_{p,g}$ on $C_\tau(V_{\eps,p})$.
On the other hand, we can consider the Dirac-class 
$[D_{p,g}]=[D_{V_{\eps,p}}]\in \KK^{\Gamma_{p,g}}(C_\tau(V_{\eps,p}),\CC)$.
It is represented by the restriction $F_{p,g}$ 
of the bounded de-Rham 
operator $F=D(1+D^2)^{-1/2}$ to $L^2(\Lambda^*_\CC(V_{\eps,p}))$.

The Kasparov product $\Theta_{p,g}\otimes_{C_\tau(V_{p,\eps})}[D_{p,g}]\in 
\KK^{\Gamma_{p,g}}(\CC,\CC)$ is represented by the pair 
$(\mathcal H_{p,g}, P_{p,g})$,
with $\mathcal H_{p,g}=L^2(\Lambda^*_\CC(V_{p,\eps}))$ and 
\begin{equation}\label{eq-prodpg}
P_{p,g}= \lambda_{\theta_{p,g}(x)}+\lambda_{\theta_{p,g}(x)}^*+
\sqrt{1-\|\theta_{p,g}(x)\|^2}F_{p,g}.
\end{equation}
(compare with the proof of Lemma \ref{lem-local} above).
Since $P_{p,g}$ is $\Gamma_{p,g}$-invariant, it restricts to an 
operator $Q_{p,g}$ on the subspace of the $\Gamma_{p,g}$-invariant vectors 
in $\H_{p,g}$. We then get well defined classes
\begin{equation}\label{kasproduct}
[(\H_{p,g}^{\Gamma_{p,g}}, Q_{p,g})]\in \KK(\C,\C)
\end{equation}
for each pair $(p,g)\in S$. We now get

\begin{proposition}\label{prop-kasproduct}
The class $\sigma_*\Deltah\hot_{\Lambda\hot \Lambdah}\alpha^*\Delta\in \KK(\C,\C)$ is equal to 
the sum 
$$\sum_{(p,g)\in \Sigma} [(H_{p,g}^{\Gamma_{p,g}}, Q_{p,g}]\in \KK(\C,\C).$$ 
\end{proposition}
\begin{proof}
Recall the operators $F\hot 1=D(1+D^2)^{-1/2}\hot 1$ and
$\theta_G$ from Definitions \ref{classdelta} and
\ref{classdeltah}. 
Under the identification,
$\mathcal{E}\hot_{\Lambda\hot \Lambdah} \ltwoformsg\cong
\oplus_{p\in P}L^2(\Lambda^*_\C F_{\eps,p})^{K_p}$ of Lemma 
\ref{lem-hilbert}, the operator $\theta_G\hot 1$ on
$\E\hot_{\Lambda\hot \Lambdah}\ltwoformsg$ corresponds to
the operator $\tilde{\Theta}$ given point-wise by
$\lambda_{\tilde{\theta}(x,g)}+\lambda_{\tilde\theta(x,g)}^*$,
with $\tilde{\theta}(x,g)=\theta(\phi(gx),x)$ and 
where $\lambda_v$ denotes exterior multiplication with $v$.
Observe
that
\begin{align*}
h\tilde{\theta}(x,g) &= d_h^x\bigl( \tilde{\theta} (h^{-1}x, \zeta (h)^{-1}gh)\bigr)
= d_h^x\bigl( \theta ( h^{-1}\phi (gx), h^{-1}x)\bigr)\\
&= \theta\bigl( \phi (gx), x\bigr)
=\tilde{\theta}(x,g),
\end{align*}
 since $\theta$ is $G$-invariant.
It follows that $\tilde{\Theta}$ 
descends to an
operator on each $L^2(\Lambda^*_{\C}F_{\epsilon,p})^{K_p}$.
Under the decomposition 
$L^2(\Lambda^*_{\C}F_{\epsilon,p})^{K_p}\cong \oplus_{g\in \Sigma_p}\H_{p,g}^{\Gamma_{p,g}}$
this operator becomes the sum 
$\oplus_{g\in \Sigma_p} (\lambda_{\theta_{p,g}(x)}+\lambda_{\theta_{p,g}(x)}^*)$ as in 
(\ref{eq-prodpg}).
Similarly, the operator  $1\hot (F\hot 1)$ descents to the sum of the de Rham operators 
$F_{p,g}$
under the decomposition $\mathcal{E}\hot_{\Lambda\hot \Lambdah}\ltwoformsg
\cong \oplus_{(p,g)\in \Sigma} \H_{p,g}^{\Gamma_{p,g}}$.
To check that the sum  of the operators
$$Q_{p,q}= \lambda_{\theta_{p,g}(x)}+\lambda_{\theta_{p,g}(x)}^*+
\sqrt{1-\|\theta_{p,g}(x)\|^2}F_{p,g}$$
on $\oplus_{(p,g)\in \Sigma}\H_{p,g}^{\Gamma_{p,g}}$ satisfies the axioms of a Kasparov product as explained in Remark \ref{rem-product}, it is enough to check that 
$T:=\oplus_{(p,g)}F_{p,g}$ is a $F\hot 1$ connection. But this follows from the description
of the isomorphism $\Phi$ of Lemma \ref{lem-hilbert}: If $\xi=\alpha\hot[h]\in C_c(X)\F_U\hot\C G$,
and if we consider $T$ as an operator on $\oplus_{p\in P}L^2(\Lambda_\C^* F_{\eps,p})^{K_p}$
via the obvious identifications, then the operator
$$F_{\xi}:= \Theta_\xi (F\hot 1)-(-1)^{\deg{\xi}\deg{F}} T\Theta_{\xi}\in 
\B\big(L^2(\Lambda_\C^*(X)\hot l^2G, \oplus_{p\in P} L^2(\Lambda_\C^* F_{\eps,p})^{K_p}\big)$$
can be described as the composition of the operator
$[\Pi(\alpha[h]\hot[h]), F_2\hot 1]\in \cK(L^2(\Lambda_\C^*(X)\hot l^2G)$ followed by a projection  
to the $L^2$-sections on a finite union of components in 
$F_{\eps}$ (which are determined by the support of $\alpha$),
and then followed by the operator $\Psi$ of Lemma \ref{lem-Hilbert}, which becomes bounded 
when restricted to 
the set of $L^2$-sections on a fixed finite number of components of $F_{\eps}$.
Thus $F_{\xi}$ is a composition of a compact operator by  bounded operators, hence it is compact.
This finishes the proof.
\end{proof}

Finally, to get the corresponding integer for the class
$\Deltah\hot_{\Lambda\hot\Lambdah}\alpha^*\Delta\in \KK(\C,\C)\cong
\Z$, we want to compute the index of the operator $H$, which is the
sum of indexes of $Q_{p,g}$ for all $(p,g)\in\Sigma$. That is,
\begin{equation}\label{indexH}
\sigma_*\Deltah\hot_{\Lambda\hot
\Lambdah}\alpha^*\Delta=\sum_{(p,g)\in\Sigma}
\operatorname{Ind}(Q_{p,g}).
\end{equation}
We do this by first computing the classes $[(\H_{p,g}, P_{p,g})]\in \oR(\Gamma_{p,g})$
and then computing from this the index of $Q_{p,g}$ as in Corollary \ref{fixedindex}.
%
 
To compute  $[(\H_{p,g}, P_{p,g})]\in \oR(\Gamma_{p,g})$, we linearize using
the exponential map so that we are considering a similar
problem in Euclidean space. So in what follow we may assume that 
$p=0$ is the origin in $\R^n$ and $V_{\eps,p}$ is some open neighborhood 
of $p=0$ in $\R^n$. By choosing $\eps$ small enough, we may further assume 
that $p=0$ is the only fixed point of the differential map $x\mapsto \phi(gx)$.
This implies that the vector field $\theta_{p,g}:V_{\eps,p}\to TV_{\eps,p}$ 
only vanishes at the point $p=0$.
The group $\Gamma_{p,g}$ acts on
$TV_{\eps,0}$ through the standard action of $\on$ on
$\R^n$. Let
$\rho_{p,g}:\Gamma_{p,g}\longrightarrow \on$ be the
corresponding representation. 

We know from Lemma \ref{lem-homotopy} that the class 
$[(\H_{p,g}, P_{p,g})]\in \oR(\Gamma_{p,g})$ only depends on the homotopy class 
of $\theta_{p,g}$, where, by Lemma \ref{lem-support} we may restrict $\theta_{p,g}$ to arbitrarily
small open balls around $0$.

Note first that under the identification of $V_{\eps,p}$ with a neighborhood 
of $0$ in $\R^n$ via the
exponential map, the metric, call it $\rho$, on
$V_{\eps,p}$ is not necessarily the flat metric coming from
$\Rn$. However the convex combination of the metric
$\rho$ and the Euclidean metric, $\rho_{\Rn}$,
gives a homotopy, $\rho_t$, between these two metrics. That, in
turn, gives a homotopy,
$\theta_t:=\left(\frac{\rho_t}{\epsilon}(d_y\rho_t)\right)(\phi(gx),x))$
of corresponding vector fields.
Therefore without loss of
generality, we may assume that the set $V_{\eps,p}$ is
equipped with the Euclidean metric and 
that $\theta_{p,g}(x)=x-(\phi\circ g)(x)$.
By calculus, $x-(\phi\circ g)(x)=\left(\Id_{\Rn}-d(\phi\circ
g)(p)\right)\cdot x+\psi(x)$ for some $\psi$ such that
$\frac{\psi(x)}{\|x\|}\rightarrow 0$ as $x\rightarrow 0$. Then,
in a small neighbourhood of $p=0$, $t\mapsto
\left(\Id_{\Rn}-d(\phi\circ g)(p)\right)\cdot x+t\psi(x)$
gives a homotopy, $\theta_{p,g}\sim \left(\Id_{\Rn}-d(\phi\circ
g)(p)\right)\cdot x=W_{p,g}X$ with $W_{p,g}=\Id_{\Rn}-d(\phi\circ g)(p)$.
It follows then from Lemma \ref{charactertheorem} together with Corollary 
\ref{fixedindex} that 
$$\Ind(Q_{p,g})=\frac{1}{|\Gamma_{p,g}|}\sum_{g\in \Gamma_{pg}}\chi_{(\rho_{p,g}, W_{p.g})}(g),$$
where $\chi_{(\rho,A)}(h)=\mathrm{sign}\,
\mathrm{det}(A|_{\operatorname{Fix}(h)})$ is the orientation character as in Definition
\ref{chi}.

Putting all together, we have the following theorem:

\begin{theorem} The pairing $\lk \widehat{[\alpha]} , \Delta\rk$ is
given by
\begin{equation}
\label{thetheorem} \lk \widehat{[\alpha]}, \Delta\rk =
\sum_{(p,g)\in\Sigma}
\frac1{|\Gamma_{p,g}|   }\sum_{h\in\Gamma_{p,g}}\chi_{(\rho_{p,g},W_{p,g})}(h).
\end{equation}
\end{theorem}

The above theorem together with the abstract Lefschetz theorem of \cite{Emerson}
proves our Lefschetz fixed point theorem, Theorem \ref{intro:mainthm}.

We now discuss an example.

\begin{example} Let $G\cong \Z\cross \Z/2\Z$
be the infinite dihedral group. It is the subgroup of
$\mathrm{Iso}(\R)$ generated by $u(x) = -x$ and $w(x) = x+1$. It has
the relation $uwu = w^{-1}$, and has two conjugacy classes of finite
subgroups $K_1
\defeq <u> = \mathrm{Stab}_G(0)$, and $K_2 \defeq <wu> =
\mathrm{Stab}_G(\frac{1}{2})$. A fundamental domain for the action
is the interval $[0,\frac{1}{2}]$. Note that $\dot{0}\not=
\dot{\frac{1}{2}} \in \dot{\R}$, where we use dot notation to
indicate orbits.

The $\K$-theory of $C_0(\R)\cross G$ is $\Z^{3}$ in dimension
$0$ and is trivial in dimension $1$. A general property of
proper actions tells us that $C(G\backslash \R) = C[0,1]$ is
strongly Morita equivalent to an ideal in $C_0(\R)\cross G$,
and one $\K$-theory generator corresponds under this strong
Morita equivalence and the inclusion of the ideal, to the class
of the unit in $C(G\backslash \R)$. We denote this class $[E]$.
The other two projections come from the $C^{*}(K_i)$, $i =
1,2$. We denote them $[p_i]$, $i = 1,2$.

Let $$\phi \colon \R \to \R, \; \;\;\; \phi (x) = -x - \frac{1}{2}.$$ Let
$\zeta\colon G \to G$ be $\zeta (u) = uw$ and $\zeta (w) = w^{-1}$. Then
$\zeta$ extends to an automorphism of $G$, and $\phi\bigl( \zeta (x)\bigr) =
g\phi (x)$ is easily checked for $g = w, u$, so that we get a covariant pair.
The map $\phi$ has one fixed orbit, which is $\dot{\frac{1}{4}}$; note that
$\phi$ itself fixes $\frac{1}{4}$. The
derivative at this point is $-1$, so that we get a positive sign attached to
this point. Since $\dot{\frac{1}{4}}$ has no isotropy in $G$, we only get a
contribution of $+1$ from this fixed orbit: the local side of the Lefschetz formula
is equal to $1$. On the global side, since $\zeta(K_1) = K_2$, there is
no tracial contribution from the summands $\Z p_1\oplus \Z p_2$, and
therefore $\mathrm{tr}_s (\alpha_*) = 1$, with $\alpha\colon C_0(\R)\cross G \to
C_0(\R)\cross G$ the induced automorphism.

For a second example, let $\zeta$ be the identity. Let $\phi$ be a
small perturbation of the identity map $\R\to \R$ which can be
roughly described as follows. Firstly, $\phi$ maps the interval
$[0,\frac{1}{2}]$ to
itself. It fixes $0$ and $\frac{1}{2}$, and has derivative zero at
both these points. It also fixes the point $\frac{1}{4}$, and has
derivative rather large at this point (in particular greater than $1$.)
Finally, $\phi$ is extended to a $G$-equivariant map $\R \to \R$ in
the obvious way.

Clearly $\phi$ is proper $G$-homotopic to the identity, so its
graded trace on $\K$-theory is $3$. It has three fixed orbits
$\dot{0}$, $\dot{\frac{1}{4}}$, and $\dot{\frac{1}{2}}$, which are
actually fixed points in $\R$. The first and third of these come
with a positive sign, and are weighted by the number of conjugacy
classes (\emph{i.e.} the number of elements) in the isotropy groups
$K_1$ and $K_2$ of these points. We thus get a contribution of
$(1+1)+(1+1) = 4$ from the first and third fixed points, and, since
$\dot{\frac{1}{4}}$ has no isotropy, and $\phi'(\frac{1}{4}) >1$, we
get a contribution of $-1$ from the second fixed point, with a net
contribution of $3$, as required.

On the other hand, if we change the above map $\phi$ just to have
now large derivatives at $0$ and $\frac1{2}$ and zero derivative at
$\frac1{4}$. Then we get a contribution $0+1$ from $0$ and also the
same from $\frac1{2}$, and $1$ from $\frac1{4}$, with a net
contribution of $3$ again.

\end{example}


\begin{thebibliography}{10}

\bibitem{AS2} M.F. Atiyah, G. Segal, \emph{The index of Elliptic Operators II},
Ann. Math., 2nd Ser., {\bf 87}, No.3.(May, 1968), 531--545.

\bibitem{AS1} M.F. Atiyah, I.M. Singer, \emph{The index of Elliptic Operators I},
Ann. Math., 2nd Ser., {\bf 87}, Issue 3 (May, 1968), 484--530.

\bibitem{Baaj-Julg} S. Baaj, P. Julg, \emph{Th\'eorie
    bivariante de Kasparov et op\'erateurs non bornes dans les
    $C^*$-modules hilbertiens}, C.R. Acad. Sci. Paris {\bf 296}
    (1983), Ser. I, 875-878.

\bibitem{BCH} P. Baum, A. Connes, N.  Higson,
    \emph{Classifying
    space for proper actions and $\K$-theory of group C*-algebras},
  Contemporary Mathematics, \textbf{167}, 241--291 (1994).

\bibitem{Blackadar} B. Blackadar, \emph{$\K$-Theory for
    Operator Algebras}, Math. Sci. Research Institute
    Publications, {\bf vol 5}, Cambridge Univ. Press, 1998.

\bibitem{EEK} S. Echterhoff, H. Emerson, H.-J. Kim.
    \emph{$\K$-theoretic  Poincar\'e duality for proper twisted
    actions.} To appear in Mathematische Annalen.

\bibitem{Em}
H. Emerson. \emph{Noncommutative Poincar\'e duality for boundary
actions of hyperbolic groups.}
 J. Reine Angew. Math. {\bf 564} (2003), 1--33.

\bibitem{Emerson}
H. Emerson.
\emph{Lefschetz numbers for $C^{*}$-algebras. }
 Preprint. arXiv:0708.4278

 \bibitem{Em1}
H. Emerson.
\emph{The Baum-Connes conjecture, noncommutative Poincar\'e duality,
and the boundary of the free group.}
 Int. J. Math. Math. Sci. {\bf 38} (2003),  2425--2445.

\bibitem{EmersonMeyer} H. Emerson, R. Meyer,
    \emph{Euler characteristics and Gysin sequences for group
    actions on boundaries.} Math. Ann., {\bf 334}, (2006), no.4,
    853-904.

\bibitem{EM} H. Emerson, R. Meyer, \emph{A descent principle
    for the Dirac dual Dirac method.} To appear in: Topology.

\bibitem{Higson-Roe} N. Higson, J. Roe, \emph{Analytic
    $\K$-Homology}, Oxford Univ. Press,(2000).

\bibitem{Kas0} G. Kasparov, \emph{Topological invariants of elliptic operators. I: K-homology},
Izv. Akad.
    Nauk SSSR Ser. Mat. {\bf 39}:4 (1975). In
    Russian; translation in Math. USSR Izvestija {\bf 9}
    (1975), 751-792.

\bibitem{Kas-operator} G. Kasparov, \emph{The operator
    $\K$-functor and extensions of $C^*$-algebras}, Izv. Akad.
    Nauk SSSR Ser. Mat. {\bf 44}:3 (1980), 571-636, 719. In
    Russian; translation in Math. USSR Izvestija {\bf 16}
    (1981), 513-572.

\bibitem{Kas1} G. Kasparov, \emph{Equivariant
      $\KK$-theory and the Novikov conjecture}, Invent. Math.
    \textbf{91},  147--201 (1988).

\bibitem{Kas2} G. Kasparov, \emph{$\K$-theory, group
    $C^*$-algebras, and higher signatures (Conspectus)}, In S.
    Ferry, A. Ranicki, and J. Rosenberg, editors, Proc. of 1993
    Oberwolfach Conference on the Novikov Conjecture, vol 226,
    LMS, 101-146, Cambridge Univ. Press, 1995.

\bibitem{KasSk} G. Kasparov and G. Skandalis. \emph{Groups
    acting properly on ``bolic'' spaces and the Novikov
    conjecture.} Annals of Math. {\bf 158} (2003), 165--206.

\bibitem{Kuc} D. Kucerovsky, \emph{The $\KK$-product of
    unbounded
    modules}, $\K$-theory {\bf 11} (1997), 17-34.
    
\bibitem{LR2} W. L\"uck annd J.M. Rosenberg,
\emph{Equivariant Euler characteristics and $K$-homology Euler classes for proper cocompact $G$-manifolds.}
 Geom. Topol. {\bf 7} (2003), 569--613.

\bibitem{Se} G. Segal, \emph{Equivariant $\K$-theory.} Publications
Math\'ematiques de l'I.H.\'E.S. {\bf 34} (1968),129--151.






\end{thebibliography}
\end{document}